\documentclass{amsart}

\usepackage{natbib}
\usepackage{epsfig}

\newcommand{\be}{\begin{equation}}
\newcommand{\ee}{\end{equation}} 
\newcommand{\abi}{_{\alpha}^{\beta_i}}

\newcommand{\kh}{\hat{\mathbf{k}}}

\newtheorem{thm}{Theorem}[subsection]

\newtheorem{defn}[thm]{Definition}
\newtheorem{remark}[thm]{Remark}

\begin{document}


\title[A Variational Free-Lagrange Method for Shallow Water]{Conservative Properties of a Variational Free-Lagrange Method for Shallow Water\footnote{Los Alamos Report LA-UR 07-7482.}}


\author[M.F. Dixon]{Matthew F. Dixon}
\address{Department of Mathematics\\
UC Davis\\
Davis, CA 95616, USA.}
\email{mfdixon@ucdavis.edu}

\author[T.D. Ringler]{Todd D. Ringler}
\address{Climate, Ocean, and Sea Ice Modeling Group\\
 Theoretical Division, Los Alamos National Laboratory\\
Los Alamos, NM 87544, USA.}
\email{ringler@lanl.gov}

\subjclass[2000]{Primary }

\keywords{Variational free-Lagrange; symplectic integrators; shallow water}

\date{}

\begin{abstract}
The variational free-Lagrange (VFL) method for shallow water is a free-Lagrange method with the additional property that it preserves the variational structure of shallow
water. The VFL method was first derived in this context by \cite{AUG84}
who discretized
Hamilton's action principle with a free-Lagrange data structure. The purpose
of this article is to assess the long-time conservation properties of the VFL method for regularized shallow water which are useful for climate simulation. Long-time regularized shallow water simulations show that the VFL method exhibits no secular drift in the (i) energy error through the application of symplectic integrators; and (ii) the potential vorticity error through the construction of discrete curl, divergence and gradient operators which satisfy semi-discrete divergence and potential vorticity conservation laws. These diagnostic semi-discrete equations augment the description of the VFL method by characterizing the evolution of its respective irrotational and solenoidal components in the Lagrangian frame. Like the continuum equations, the former exhibits a $\text{div}^2\mathbf{U}$ term which indicates that the flow has a very strong tendency towards a purely rotational state. 

Numerical results
show (i) the preservation of shape and strength of an initially radially symmetric vortex pair in purely rotational regularized shallow water and (ii) how the Voronoi diagram retains the history of the flow field and (iii) that energy is conserved to $\mathcal{O}(\Delta^2)$ and potential vorticity error to within $5\%$ with no secular growth over a 50 year period. 
\end{abstract}

\maketitle                   

\section{Introduction}
\paragraph{\emph{Motivation}} The derivation and analysis of conservative numerical methods for global
climate modeling is a challenging and active area of research. The traditional
approach to deriving such methods requires the numerical analyst or engineer to exercise some level of their own bias in choosing which conservation laws to enforce in the discrete model. To remedy this and systematically derive the numerical method
from a conferred mathematical understanding, \cite{SALMON83} pursues the \emph{variational
approach} for geophysical modeling. The power of the variational approach rests upon Noether's theorem
which states that each continuous symmetry of Hamilton's action principle
exhibits a corresponding conservation law.  Consequently, the numerical analyst
need only choose a suitable discretization of Hamilton's action principle,
with the confidence that not only will the resulting semi-discrete equations of motion preserve
geometric structure, manifesting in energy conservation, but that these equations
will conserve additional quantities  corresponding to the symmetries of the
discrete action principle.

\paragraph{\emph{The free-Lagrange method}} \cite{AUG84} presents a variational extension of the free-Lagrange
method for rotating shallow water- a method which by itself exhibits a number of favorable properties. In contrast with the majority  of methods for shallow water, which use a fixed grid for transport computations,  the free-Lagrange method builds the transport into the evolution of particle positions. This step avoids introducing diffusion into the transport- an artefact of mixing the transport computed from fluxes between cells with the remaining fluid. Consequently, the free-Lagrange method exhibits local transport conservation laws such as fluid mass and tracer concentration.

Secondly, the free-Lagrange method represents the layer thickness over a \emph{Voronoi diagram}. This approach is attractive from both a theoretical and practical perspective. \cite{SERRANO05}
shows that use of a Voronoi diagram in Hamilton's discrete action principle permits translational
and rotational symmetries and hence respective linear and angular momentum
conservation laws for the discrete Euler equations. The Voronoi diagram avoids \textit{mesh-tangling} problems
which seriously compromise the performance of fixed connectivity mesh based Lagrangian methods. This does not mean that use of a Voronoi diagram is not without its own computational challenges. The generation of an entirely new Voronoi diagram at each time step is computationally prohibitive in three dimensions. Also, computational complexity manifests from the absence of constraints on the mesh connectivity.

Some of these computational efficiency shortcomings have been partially addressed by alternative approaches cited in the literature. \cite{HARLEN95} present a split Eulerian-Lagrangian scheme for viscoelastic fluids that retains the nodes as material points and reconnects them to produce the optimal Delaunay triangulation. This offers an advantage over the free-Lagrange method in that the discrete fluid equations can then be solved using Galerkin finite element approximations. \cite{PETERA96} presents a Galerkin-Lagrange finite element approximation of a shallow water model for tidal flow which attempt to resolve fluctuations between dry and wet regions. Their basis for node adjustment is a physical one. The nodes of the mesh are only moved if the neighbourhood of the node changes state between wetted and dried. This approach avoids excessive mesh distortions but requires a smoothing operator to eliminate numerical oscillations introduced by succesive modifications of the mesh in very shallow regions.

\paragraph{\emph{Variational methods}} Despite their computational efficiency, neither of these two previously described approaches preserve variational structure- a distinguishing feature of the variational free-Lagrange (VFL) method derived by \cite{AUG84}. We trace the development of variational numerical methods in geophysical modeling back to \cite{BUNEMAN82}, who describes a two-dimensional Eulerian model based upon the Clebsch representation of Hamilton's principle. He advocates a method whose stability appears to be contingent on the invariance of the discrete Hamilton's principle under particle relabelling. It is now well known, however, that the discrete Lagrangian is not invariant under continuous transformations of the labels - there is no particle relabelling symmetry \citep[see][]{RIPA81,SALMON82}. The symmetry gives Ertel's theorem of conservation of potential vorticity which in turn recovers the known connection with Kelvin's circulation theorem along surfaces of constant entropy. In the presentation of a numerical analogue of a variational shallow water \emph{blob} model, \cite{SALMON83} conjectures that potential vorticity conservation might still be conserved through modifying the discrete variational principle. Specifically he suggests constraining the angular component of the particle velocities. This approach is difficult to efficiently implement, however, and remains unsubstaniated.  

More computationally feasible variational approaches have emerged from classical Lagrangian numerical methods. Most notably, variational formulations of smooth particle hydrodynamics (SPH) methods have been developed \citep[see][and the references by the same first author therein]{BONET05}. In this work, the authors actually make use of the variational principle to extend the classical SPH method by deriving spatially dependent smoothing lengths (a limiting feature in conventional SPH methods is the use of a fixed smoothing length). There appears to be a notable advantage in using a variational SPH for highly compressible fluids.


\paragraph{\emph{Symplectic integrators}} The literature that we have surveyed on the VFL method does not consider preservation of the symplectic structure under time discretization of the particle Euler-Lagrange equations. Our primary contribution is to exploit
the symplectic structure preserved by the dispersively regularized shallow water equations. Dispersive regularization avoids loss of symplecticity of symplectic time-steppers \citep[see][for a review of symplectic
methods for particle methods]{DIXON04} by slowing high frequency error modes which would otherwise cause numerical instabilities. We demonstrate by numerical experiment
that the combination of the regularized VFL method with
an explicit symplectic time-stepper exhibits
no secular drift in the energy error for moderate time steps and coarse spatial resolutions. This property is desirable for long-time geophysical simulations which rely on moderate time step sizes and coarse spatial resolutions for computational tractability.

\paragraph{\emph{Overview}} We shall begin by revisiting the free-Lagrange method and subsequently deriving the semi-discrete VFL equations for regularized shallow water in Section \ref{sect:vfl_2d}. In Section \ref{sect:SYMSTT} we present the main feature of the paper, an explicit symplectic integrator for preserving the symplectic structure of the semi-discrete VFL equations for regularized shallow water. The form of the regularization operator and its implementation is described in Section \ref{sect:reg}. Using the semi-discrete VFL equations, we derive the semi-discrete vorticity and divergence equations in Sections \ref{sect:vor} and \ref{sect:div} and show that the potential vorticity is only conserved if the discrete curl operator is judiciously chosen so that its operation on the layer thickness gradient vanishes. In Section \ref{sect:num_exp}, we present numerical results which show the scaling laws of potential vorticity and energy. The VFL method for 1D rotating shallow water is shown to conserve energy to the order of the integrator and exhibit the geostrophic adjustment mechanism of rotating
shallow water. Numerical simulations illustrate the motion of a vortex pair in 2D rotating shallow water and the configuration of the Voronoi diagram. We also show the corresponding energy and potential vorticity conservation properties over a 50 year period. Section \ref{sect:sum} concludes.

\section{The Free-Lagrange Method}
Consider a particle representation of a fluid at time $t$ in which $N^2$ particles
or \emph{sites} $\mathbf{X}(t)=\{\mathbf{X}_1(t),\dots,\mathbf{X}_{N^2}(t)\}$ are
individually labelled by $\alpha$. Each particle is inside a 
polygon of area $A_{\alpha}(t)$ containing the 
set of distinct points $\mathbf{X}_{\beta}(t)$ closer to $\mathbf{X}_{\alpha}(t)$ than to any other site. The polygon is referred to as a \emph{Voronoi cell}
and the set of all these closed cells is the \emph{Voronoi diagram}. A hexagonal Voronoi cell is shown in Figure 1 together with the
specification
of a local index for referring to neighboring particles. Note that the assumption that the sites are distinct does not prohibit the formation of shocks - characteristics may still collide as they are not the particle trajectories \citep[see, for example,][ for an explanation of shock formation]{WHITHAM74}. 

\begin{figure} [!ht]
\begin{center}
\includegraphics[width=0.7\linewidth,height=65mm]{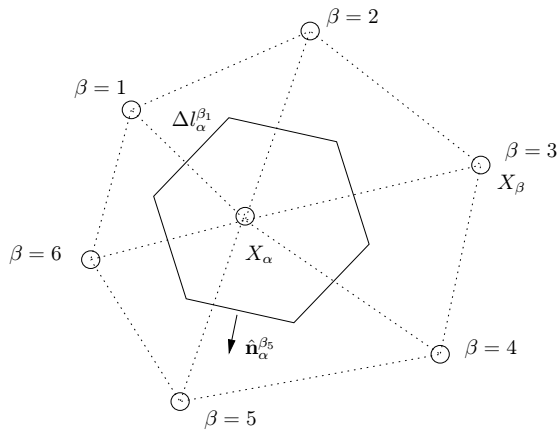}
\end{center}
\caption{\scriptsize{A hexagonal ($n_e=6$) Voronoi cell (outlined with the solid line) containing the particle
with label $\alpha$. Each cell edge is of length $\Delta l\abi$ and
indexed by $\beta_i$. The surrounding Delaunay triangulation, dual to the Voronoi diagram, is shown with the dotted line.}}
\label{fig:1}
\end{figure}

The \emph{free-Lagrange method} simply uses the particles to represent the material velocity field $\mathbf{U}_{\alpha}$ and the Voronoi diagram to construct the piecewise
constant cellwise layer thickness $\bar{h}^{n}_{\alpha}:=\bar{h}(\mathbf{X}_{\alpha},t^n)$ which, in discrete time, is computed from the cell mass conservation law

\be \label{eq:mass_cons}
\bar{h}^{n+1}_{\alpha}=\frac{A^n_{\alpha}}{A^{n+1}_{\alpha}}\bar{h}^n_{\alpha}.
\ee

\section{A Variational Free-Lagrange method for 2D Shallow Water}\label{sect:vfl_2d}

Discretization of Hamilton's
action principle for rotating shallow water with this free-Lagrange fluid representation is the critical step in deriving a VFL method. The semi-discrete material description of Hamilton's action principle for
2D rotating shallow water in a f-plane is expressed in terms of the particle material
velocities and piecewise constant cellwise layer thickness.

\begin{defn}[The Semi-Discrete Hamilton's Action for Regularized Shallow Water]
The semi-discrete action principle for rotating shallow water over bottom topography in a f-plane is
\be \label{eqn:disc_ap}
\begin{split}
\mathcal{S}_d &=\frac{1}{2}\int_{t_a}^{t_b} dt\sum_{\alpha} m_{\alpha}\left(|\mathbf{U}_{\alpha}(t)|^2 + 2\mathbf{R}_{\alpha}\cdot\mathbf{U}_{\alpha}\right)
- V,\qquad V=\frac{g}{2}\sum_{\alpha} m_{\alpha}\left(\tilde{h}(\mathbf{X}_{\alpha},t)+2\bar{b}_{\alpha}\right),
\end{split}
\ee
with fixed end points $\delta X_{\alpha}(t_a)=\delta X_{\alpha}(t_b)=0,~~\forall \alpha$. $\tilde{h}$ is a dispersive regularization of the layer thickness\footnote{The details of the dispersive regularization operator are given in Section \ref{sect:reg}.}.
\end{defn}

Under the f-plane approximation, $\nabla\times \mathbf{R}=f_0\kh$,~$\kh$ is the unit vector in
the opposite direction of gravity and $f_0$ is the Coriolis parameter which is given by $f_0=2|\omega|$, where $\omega$ is the angular velocity of rotation relative
to an inertial frame. $g$ is the gravitation constant, $\bar{b}_{\alpha}$ is
the bottom topography, assumed to be piecewise constant over cell $\alpha$
and $V$ is the potential energy.

Stationarity of the discrete action principle is equivalent to the semi-discrete Euler-Lagrange particle equations
\be \label{eqn:PIC_EL}
\dot{\mathbf{U}}_{\alpha}=g\sum_i\left[\tilde{h}\right]\abi\mathbf{dn}\abi - f_0\kh\times\mathbf{U}_{\alpha},\qquad
 \dot{\mathbf{X}}_{\alpha}=\mathbf{U}_{\alpha},
\ee
where $\mathbf{dn}\abi:=\hat{\mathbf{n}}\abi\Delta l\abi$ is the outward normal vector at edge $\beta_i$ of cell $\alpha$, $\Delta l\abi$ is the edge length and the averaging operator $\left[\cdot\right]\abi=\frac{1}{2m_{\alpha}}\left(m_{\alpha}(\cdot_{\alpha})~+~ m_{\beta_i}(\cdot_{\beta_i})\right)$ evaluates a cellwise scalar quantity over this edge
as the mass weighted mean of that quantity over cell $\alpha$ and cell $\beta_i$. These weights are constants since, by construction, the mass is conserved in each cell. The right hand side of \eqref{eqn:PIC_EL} gives the form of the cellwise discrete gradient operator 

\be \label{eqn:DISC_GRAD}
\text{grad}\left(\tilde{h}(\mathbf{X}_{\alpha},t)\right):=-\sum_i\left[\tilde{h}\right]\abi\mathbf{dn}\abi,
\ee
which is a discrete analogue of the weak form of the continuum gradient operator uniquely defined over the interior of a closed cell whose area tends to zero. This discrete gradient operator is given in weak form, \eqref{eqn:PIC_EL} applies pointwise to the material velocity vector field which is constructed at particle positions only.

\begin{remark}
We note that the discrete gradient operator \eqref{eqn:DISC_GRAD} does not define a gradient along the cell edges. The prognostic discrete divergence and potential vorticity conservation laws corresponding to the weak form of continuum conservation laws, given in Section \ref{sect:pv}, require a definition of the gradient of the regularized layer depth along cell edges. We look instead to the linear interpolant of the layer depth constructed over the Delaunay triangulation, dual to the Voronoi diagram, to evaluate a gradient here for prognostic purposes only. We emphasize that the evaluation of the gradient along cell edges just enables a more complete characterization of the conservative properties of the VFL method and is not intrinsic to the VFL method.
\end{remark}

\section{A Symplectic Time Stepping Scheme} \label{sect:SYMSTT}

The semi-discrete Euler-Lagrange particle equations \eqref{eqn:PIC_EL} are symplectic and
conserve energy. A symplectic time-stepper, such as the explicit
second order St\"ormer-Verlet scheme, in particle position and canonical momentum co-ordinates $(\mathbf{X}^n_{\alpha},
\mathbf{P}^n_{\alpha}=m_{\alpha} \mathbf{U}^n_{\alpha})$ at
time $t^n$ 
\begin{equation}\label{eqn:HAM}
\begin{split}
\mathbf{P}_{\alpha}^{n+\frac{1}{2}}&=\mathbf{P}_{\alpha}^n-\frac{\Delta t}{2}V_{\mathbf{X}}(\mathbf{X}_{\alpha}^n),\\
\mathbf{X}_{\alpha}^{n+1}&=\mathbf{X}^n_{\alpha}+\Delta t \frac{\mathbf{P}_{\alpha}^{n+\frac{1}{2}}}{m_{\alpha}},\\
\mathbf{P}_{\alpha}^{n+1}&=\mathbf{P}_{\alpha}^{n+\frac{1}{2}}-\frac{\Delta t}{2}V_{\mathbf{X}}(\mathbf{X}_{\alpha}^{n+1}),
\end{split}
\end{equation}
\noindent preserves the symplectic two-form $\omega^n=\sum_{\alpha=1}^{N^2} d\mathbf{X}^n_{\alpha}\wedge d\mathbf{P}^n_{\alpha}$, where
$V_{\mathbf{X}}$ denotes the gradient of the potential energy.

\vspace{.5cm}
\paragraph{\emph{Discussion}} Integrators which conserve the Hamiltonian are often regarded as a strong manifestation of unconditional numerical stability \citep{LEWIS94} and therefore permit the use of large time steps. Symplectic integrators, on the other hand, conserve a modified Hamiltonian which is derived by backward error analysis.  When used as time steppers for variational free-Lagrange methods, both the particle resolution and the time step size must be chosen to ensure proximity of the modified Hamiltonian to the exact Hamiltonian. Otherwise, particles collide resulting in singularities in the Jacobian of the surrounding Voronoi cells. In such cases, the particles must be remapped and the variational structure is lost. 

Dispersive regularization of the layer depth prevents particles colliding with a time step size generally no smaller than required by semi-implicit semi-Lagrangian methods. These methods, despite being dissipative, are preferred by practitioners as they permit larger time steps than Eulerian methods. The reader is referred to \cite{REICH06} for details of dispersive regularization methods for particle mesh approximations of shallow water and their comparison to the dispersive properties of semi-implicit time stepping.

The principal advantage of using the variational free-Lagrange method combined with symplectic integrators over non-symplectic time steppers for semi-Lagrangian methods is the general absence of secular drift in the energy error. The ever diversifying set of computational models which crucially rely on the ability to tractably simulate long-time dynamics without secular numerical dissipation motivates the development of new computationally competitive geometric numerical methods such as the VFL method. Global climate models are arguably one of the strongest candidates in this regard- implementers are challenged by the need to tractably resolve complex geophysical dynamics over long periods.

\section{Dispersive Regularization} \label{sect:reg}

The modification of the semi-discrete action principle, to include a dispersively regularized potential, ensures that the corresponding time stepping scheme \eqref{eqn:HAM} remains stable for moderately large time steps at computationally tractable resolutions for long-time global climate simulations. This is frequently accomplished by approximating the inverse of a Helmholtz elliptic operator 

\be \label{eqn:IHH}
\tilde{h}=(1-\hat{\alpha}^2\nabla)^{-1}h,
\ee
which includes a characteristic dispersive regularization length scale $\hat{\alpha}$. We construct a conforming finite element approximation of this operator on the Delaunay triangulation $\mathcal{T}_d$ which is dual to the Voronoi diagram over the domain $\Omega$.

More precisely, find $\tilde{h}_d\in V(\tilde{h}_d)$ with $a(\tilde{h}_d,\phi)=b(h_d, \phi)$ for all $\phi$ in the space of all linear functions. $a(,)$ and $b(,)$
are bilinear forms respectively defined on the finite dimensional vector spaces $V(\tilde{h}_d)$ and $V(h_d)$ where
\begin{equation}
V(h_d):=\{h_d\in C(\bar{\Omega}) ~|~ h |_{\mathcal{K}}\in \mathcal{P}_1(\mathcal{K})
~\forall \mathcal{K}\in \mathcal{T}_d\subseteq H^1_0(\Omega)\}.
\end{equation}
Let $\tilde{p}_{\mathcal{K}},~p_{\mathcal{K}}\in\mathcal{P}_1$ be the
polynomial interpolating functions on triangle $\mathcal{K}$ where

\begin{equation}
\tilde{p}_{\mathcal{K}}(x,y)=\sum_i \tilde{h}(x,y)N_i(x,y)\qquad\text{and}~~~p_{\mathcal{K}}(x,y)=\sum_i h(x,y)N_i(x,y),
\end{equation}
and the locally indexed basis functions over triangle $\mathcal{K}$ take the form
\begin{equation}
\begin{split}
N_1(x,y)&=[(x_2y_3-x_3y_2) + (y_2-y_3)x + (x_3-x_2)y)]/2A,\\
N_2(x,y)&=[(x_3y_1-x_1y_3) + (y_3-y_1)x + (x_1-x_3)y)]/2A,\\
N_3(x,y)&=[(x_1y_2-x_2y_1) + (y_1-y_2)x + (x_2-x_1)y)]/2A,\\
\end{split}
\end{equation}
where $A$ denotes the area of triangle $\mathcal{K}$. Defining $\tilde{h}_d(x,y)=\tilde{p}_{\mathcal{K}}(x,y),~h_d(x,y)=p_{\mathcal{K}}(x,y)~~~\forall x \in \mathcal{K}\in\mathcal{T}_d$ and constructing the sparse stiffness matrix $A$ and mass matrix $M$ representing $a(\tilde{h}_d,N_i)=b(h_d,N_i),~\forall
i$ to give the sparse linear system
\be\label{eqn:linear_iso}
A\tilde{h}=M h,
\ee
where the matrices are defined from the weak form of \eqref{eqn:IHH}

\begin{equation}
\begin{split}
A_{ij}=a(N_i,N_j)&=\sum_{\mathcal{K}\in\mathcal{T}_d} \int_{\mathcal{K}}
d^2x N_i N_j + \hat{\alpha}^2(N_i)_{,k}(N_j)_{,k},\\
M_{ij}=b(N_i,N_j)&=\sum_{\mathcal{K}\in\mathcal{T}_d} \int_{\mathcal{K}}
d^2x N_iN_j.
\end{split}
\end{equation}

\paragraph{Computational procedure} At each time step, the linear system \eqref{eqn:linear_iso} is solved by inverting the sparse stiffness matrix to give the regularized layer thickness at each Voronoi site. This layer thickness is a prognostic variable used to evaluate the discrete cellwise gradient operator \eqref{eqn:DISC_GRAD} in the right hand side of \eqref{eqn:PIC_EL}. We outline the basic mechanism for applying the regularization to the VFL method: 

\vspace{1.0cm}
\begin{enumerate}
\item The unregularized layer thickness is updated at each time step according to the mass conservation law \eqref{eq:mass_cons};\\
\item The linear interpolant of the updated unregularized layer thickness $h_d$ is constructed over the Delaunay triangulation; \\
\item The Helmholtz linear system \eqref{eqn:linear_iso} is solved for the linear interpolant of the regularized layer thickness; and\\
\item The discrete gradient \eqref{eqn:DISC_GRAD} is evaluated from the nodal values of the regularized layer thickness.
\end{enumerate}

The linear interpolant of the layer thickness over the triangulation uniquely defines the gradient of the regularized layer thickness at the edge between neighboring sites with global indices $\alpha$ and $\beta_i$

\be\label{eqn:layer_thickness}
\left(\text{grad}\left(\tilde{h}_{\alpha}\right)\right)\abi:=\frac{\hat{\mathbf{n}}\abi}{d\abi}\left(\tilde{h}_{\beta_i} -\tilde{h}_{\alpha}\right),
\ee
where $d\abi$ is the Euclidean distance between sites $\alpha$ and $\beta_i$.  This form of the gradient operator will now be used to derive the diagnostic semi-discrete potential vorticity and divergence equations.

\section{The Shallow Water Potential Vorticity and Divergence Equations}\label{sect:pv}

Geophysical flows are generally dominated by potential vorticity dynamics, a property which together with the divergence underlies the prognostic momentum equation. We seek conditions under which the discrete momentum equations have underlying potential vorticity and divergence dynamics which are consistent with Ertel's theory of potential vorticity conservation and the continuum divergence equation respectively. To this end we derive the semi-discrete potential vorticity and divergence equations from the discrete momentum equations and identify spurious terms arising from the semi-discretization.

The curl and divergence of velocity respectively kinematically equate to the solenoidal and irrotational components of the Helmholtz decomposed velocity field. For this reason, the curl and divergence evolution equations are typically derived together and describe how each corresponding component of the velocity field evolves. Moreover, the decomposition has the added benefit of largely separating the temporally fast modes characteristic of the irrotational flow from the temporally slow modes characteristic of the solenoidal flow. 

The formulation of the semi-discrete potential vorticity and divergence equations depends on the definition of discrete curl, div and grad operators. We have already shown that semi-discretization of Hamilton's action principle with the free-Lagrange data structure defines the discrete grad operator over the interior of the Voronoi cells. The remaining discrete operators, including the definition of the layer thickness gradient along cell edges is not defined by the VFL method. 

Following \cite{RANDALL02}, we define the 
normal component of the discrete curl operator by analogy with the derivation of the discrete divergence operator given by e\eqref{eqn:chpt_4_div}. The Voronoi diagram is assumed to rest in the X-Y plane so that the vertical unit vector $\kh$ is normal to the plane. From Stoke's theorem we arrive at

\be
\frac{1}{A_{\alpha}}\int_{A_{\alpha}}\kh\cdot\nabla\times\mathbf{U}dA=\frac{1}{A_{\alpha}}\oint_{\partial
A_{\alpha}}\hat{\mathbf{\tau}}\cdot \mathbf{U} dl=\frac{1}{A_{\alpha}}\sum_{i=1}^{n_e}\hat{\mathbf{\tau}}\abi\cdot\oint_{\partial
A\abi}\mathbf{U}dl,
\ee
which approximates to

\be \label{eqn:chpt_4_div}
\frac{1}{A_{\alpha}}\sum_{i=1}^{n_e}\hat{\mathbf{\tau}}\abi\cdot\int_{\partial
A\abi}\mathbf{U}dl\approx \frac{1}{A_{\alpha}}\sum_{i=1}^{n_e}\mathbf{d\tau}\abi\cdot(\mathbf{U})\abi:=\kh\cdot\text{curl}(\mathbf{U}_{\alpha}),
\ee
where $\mathbf{d\tau}\abi:=\hat{\mathbf{\tau}}\abi\Delta l\abi$. The discrete curl operator applies to a cell edgewise vector field and not to the vector field constructed at particle positions. As previously mentioned, because the discrete gradient operator \eqref{eqn:DISC_GRAD} is not uniquely defined over cell edges, we instead evaluate the gradient of the layer thickness at the cell edges using the interpolant constructed over the dual mesh.

The additional discrete operators needed to establish discrete divergence and potential vorticity laws are the discrete div and curl operators which respectively take the form

\be
\begin{split}
\text{div}\left(\cdot_{\alpha}\right)&:=\frac{1}{A_{\alpha}}\sum_{i=1}^{n_e} \mathbf{dn}\abi\cdot\left(\cdot\abi\right),\\
\kh\cdot \text{curl}\left(\cdot_{\alpha}\right)&:=\frac{1}{A_{\alpha}}\sum_{i=1}^{n_e} \mathbf{d\tau}\abi\cdot\left(\cdot\abi\right),\\
\end{split}
\ee
and are defined in cell $\alpha$ by evaluating the argument along the cell edges.

\subsection{The semi-discrete potential vorticity equation} \label{sect:vor}

The vertical component of the discrete curl of the semi-discrete Euler-Lagrange equation
over cell $\alpha$ takes the form
\be \label{eqn:curl_pme}
\begin{split}
\kh\cdot \text{curl}(\dot{\mathbf{U}}_{\alpha})&=-g\kh\cdot \text{curl}\left(\text{grad}(\tilde{h}(\mathbf{X}_{\alpha},t))\right) - f_0\kh\cdot \text{curl}(\kh\times \mathbf{U}_{\alpha}).
\end{split}
\ee
Using the definition of \emph{curl}, after several calculations provided in Section \ref{sect:vor_calc}, the left hand side of \eqref{eqn:curl_pme} becomes
\be \label{eqn:div_calc}
\begin{split}
&\kh\cdot \text{curl}(\dot{\mathbf{U}}_{\alpha})
=\frac{D}{Dt}\left(\kh\cdot \text{curl}(\mathbf{U}_{\alpha})\right) + \text{div}(\mathbf{U}_{\alpha})
\kh\cdot\text{curl}(\mathbf{U}_{\alpha}),\\
\end{split}
\ee
where 
\be
\frac{D}{Dt}\mathbf{d\tau}\abi=\frac{D}{Dt}\left(\hat{\mathbf{\tau}}\abi\Delta
l\abi\right)=\Delta\frac{D}{Dt}(\hat{\mathbf{\tau}}l)\abi=\Delta U\abi:=U_{\alpha}^{\beta_{i+1}}-U_{\alpha}^{\beta_{i-1}}.
\ee 
Substituting the definition of the relative vorticity
$\zeta_{\alpha}=\kh\cdot\text{curl}(\mathbf{U}_{\alpha})$ into the last line of \eqref{eqn:div_calc}, we arrive at 
\be
\kh\cdot \text{curl}\frac{D \mathbf{U}_{\alpha}}{Dt}=\frac{D}{Dt}\zeta_{\alpha} + \zeta_{\alpha}\text{div}(\mathbf{U}_{\alpha}).
\ee
Taking this expression and substituting it into \eqref{eqn:curl_pme} gives the semi-discrete vorticity equation
\be \label{eqn:semi_vort}
\frac{D}{Dt}\zeta_{\alpha} - \zeta_{\alpha}\text{div}(\mathbf{U}_{\alpha}) + f_0\kh \cdot \text{curl}\left(\mathbf{U}^{\perp}\right)=0,\\
\ee
\noindent where the right-hand side vanishes because the discrete \emph{curl} of the edgewise gradient \eqref{eqn:layer_thickness} evaluates to 
\be \label{eq:curl_grad}
\begin{split}
\kh\cdot\text{curl}(\text{grad}(\tilde{h}_{\alpha})) &=
\frac{1}{A_{\alpha}}\sum_i
\left(\text{grad}(\tilde{h}_{\alpha})\right)\abi\cdot\mathbf{d\tau}\abi\\
&=\frac{1}{A_{\alpha}}\sum_i
\frac{\hat{\mathbf{n}}\abi}{d\abi}\left(\tilde{h}_{\beta_i} -\tilde{h}_{\alpha}\right)\cdot\mathbf{d\tau}\abi\\
&=0.
\end{split}
\ee
Using the following properties of the discrete curl and div operators
\be
\begin{split}
\kh\cdot\text{curl}(\mathbf{U}_{\alpha}^{\perp})&=
\frac{1}{A_{\alpha}}\sum_{i=1}^{n_e} \mathbf{d\tau}\abi\cdot\left(\mathbf{U}_{\alpha}^{\perp}\right):=
\frac{1}{A_{\alpha}}\sum_{i=1}^{n_e} \mathbf{d\tau}\abi\cdot\left(\mathbf{U}_{\alpha}\times\hat{\mathbf{k}}\right)\\
&=\frac{1}{A_{\alpha}}\sum_{i=1}^{n_e} \mathbf{dn}\abi\cdot\left(\mathbf{U}_{\alpha}\right)
=\text{div}(\mathbf{U}_{\alpha}),
\end{split}
\ee
the reader can easily confirm that \eqref{eqn:semi_vort} simplifies to

\be \label{eqn:sd_vor}
\frac{D}{Dt}\left(\zeta_{\alpha} + f_0\right) + \left(\zeta_{\alpha} + f_0\right)\text{div}(\mathbf{U}_{\alpha})=0.
\ee

\paragraph{\emph{Potential vorticity}} The semi-discrete potential vorticity equation follows immediately from the semi-discrete vorticity equation and the semi-discrete mass continuity equation

\be
\frac{D}{Dt}\bar{h}_{\alpha}=-\bar{h}_{\alpha}\text{div}(\mathbf{U}_{\alpha}).
\ee
Substituting this expression for the discrete divergence of the particle velocity $\mathbf{U}_{\alpha}$ into \eqref{eqn:sd_vor} gives the semi-discrete potential vorticity conservation law 

\be
\frac{D}{Dt}\left(\zeta_{\alpha} + f_0\right) + \left(\frac{\zeta_{\alpha} + f_0}{\bar{h}_{\alpha}}\right)\frac{D}{Dt}\bar{h}_{\alpha}=\frac{D}{Dt}\left(\frac{\zeta_{\alpha} + f_0}{\bar{h}_{\alpha}}\right)=0.
\ee

We observe that this semi-discrete potential vorticity conservation law is conditional upon the discrete curl and gradient operators satisfying \eqref{eq:curl_grad}.

\subsection{The semi-discrete divergence equations} \label{sect:div}

We now repeat the previous steps for deriving the semi-discrete potential
vorticity equation, only this time, taking the discrete divergence of the
semi-discrete Euler-Lagrange equations. The discrete divergence of the semi-discrete Euler-Lagrange equation
over cell $\alpha$ takes the form
\be \label{eqn:div_pme}
\begin{split}
\text{div}(\dot{\mathbf{U}}_{\alpha})&=-g\text{div}\left(\text{grad}(\tilde{h}(\mathbf{X}_{\alpha},t))\right) - f_0\text{div}(\kh\times \mathbf{U}_{\alpha}).
\end{split}
\ee
Using the definition of \emph{div}, after several calculations provided in \ref{sect:div_calc}, the left hand side of \eqref{eqn:div_pme} becomes
\be \label{eqn:chpt_div_calc}
\begin{split}
&\text{div}(\dot{\mathbf{U}}_{\alpha})=\frac{1}{A_{\alpha}}\sum_i \dot{\mathbf{U}}\abi\cdot \mathbf{dn}\abi=\frac{D}{Dt}\text{div}(\mathbf{U}_{\alpha}) + \text{div}(\mathbf{U}_{\alpha})^2
- \Gamma_{\alpha}.\\
\end{split}
\ee
Denoting the divergence of $\mathbf{U}_{\alpha}$ as $\delta_{\alpha}:=\text{div}(\mathbf{U}_{\alpha})$, \eqref{eqn:chpt_div_calc} becomes 
\be \label{eqn:semi-discrete_div}
\text{div}(\frac{D}{Dt}\mathbf{U}_{\alpha})=\frac{D}{Dt}\delta_{\alpha} + \delta_{\alpha}^2 - \Gamma_{\alpha}.\\
\ee

\noindent So the semi-discrete divergence equation over cell $\alpha$ is 
\be \label{eqn:se_div}
\frac{D}{Dt}\delta_{\alpha} + \delta_{\alpha}^2 - \Gamma_{\alpha} + g\text{div}(\text{grad}(\tilde{h}_{\alpha})) - f_0\text{div}(\mathbf{U}^{\perp}_{\alpha})=0.\\
\ee

\begin{remark}
We compare the semi-discrete shallow water divergence term in \eqref{eqn:semi-discrete_div} with its continuous
form given by

\begin{equation}
\begin{split}
\text{div}\left(\frac{D}{Dt}\mathbf{U}\right)&=\partial_t\delta
+ \text{div}(\mathbf{U}\cdot\nabla\mathbf{U})\\
&=\partial_t\delta
+ \partial_X(U\partial_X U + V\partial_Y U) + \partial_Y(U\partial_X V + V\partial_Y
V)\\
&=\partial_t\delta
+ (\partial_X U)^2 + (\partial_Y V)^2 + 2\partial_Y U\partial_X V + U\partial_X^2U
+ V\partial_Y^2V + V\partial_{XY}U + U\partial_{XY}V\\
&=\partial_t\delta
 + (\partial_X U)^2 + (\partial_Y V)^2 + 2\partial_Y U\partial_X V + \mathbf{U}\cdot\nabla\delta\\
 &=\frac{D}{Dt}\delta
+ \delta^2 +2(\partial_XV\partial_YU-\partial_XU\partial_YV)\\
 &=\frac{D}{Dt}\delta
+ \delta^2 - 2\text{det}(\nabla \mathbf{U}).
\end{split}
\end{equation}

Intruiguingly, because the Lagrangian form of the continuum shallow water divergence equation is less visited by the literature, we find the need to emphasize the physical significance of the $\text{div}^2\mathbf{U}$ term, which has a very strong tendency to force the fluid back to a purely rotational state.  The discrete variant of this term also appears in the semi-discrete divergence equations. Through appealing to further analysis of this term's effect on continuum vortex dynamics can we formulate an analogous statement in support of correctly resolved vortex dynamics in simulation.  We also note that $\Gamma_{\alpha}$ only evaluates to $\frac{1}{A_{\alpha}}\oint_{\partial
A_{\alpha}}
\mathbf{U}\cdot d(\mathbf{U}^{\perp})\equiv 2\text{det}(\nabla \mathbf{U})$ in the continuum limit.  So like the semi-discrete potential vorticity, we are able to identify analogous physically significant terms, but the extent to which these quantities are conserved in simulation is governed by the consistency of the choosen discrete operators with the approximation solution fields.
\end{remark}

\section{Numerical Experiments}\label{sect:num_exp}

This Section describes three numerical experiments using the VFL to simulate the rotating regularized shallow water equations on a periodic domain. The purpose of the first experiment
is to investigate the conservative properties of the VFL method for rotating 1D
shallow water. In the second experiment, we consider the problem of geostrophic
adjustment of rotating 1D shallow water and show that the VFL method produces results
which are consistent with the theory of geostrophic adjustment established
by \cite{ROSSBY38}. 
This theory describes the physical mechanism by which
perturbed rotating shallow water recovers to a geostrophically balanced state.
Layer depth perturbations $h'$ of the scale $L'<<L_D$, where $L_D=\sqrt{gH_0}/f_0$
is the Rossby deformation radius for shallow water in a f-plane, result in
gravity waves which propagate energy and momentum away from the source leaving behind a geostrophically balanced flow. It is precisely because this mechanism is consistent with the advection law of potential vorticity that geophysical modelers use this experiment to assess the integrity of the discrete shallow water velocity and layer depth equations. These two experiments are initialized by perturbing the Voronoi cell masses with a Gaussian perturbation and the velocity field is initially zero. 

In the third experiment, we simulate the motion of a vortex pair in purely rotational 2D shallow water flow and monitor the energy and potential vorticity conservative properties over a much longer period of 50 years. We verify that the geostrophic balances state is preserved - the vortex pair undergoes pure rotation without shape deformation or a change in strength of the poles. Over the course of the simulation, the Voronoi diagram appears to retain the history of the flow field and deforms in the wake of the vortex pair. We also show that the relative energy and potential vorticity error does not exhibit secular drift. 

Before describing each of the experiments in more detail, we mention a few pertinent details common to all experiments.  The smoothing parameter $\hat{\alpha}$ is choosen to be at least equal to the reference grid scale $dx$- that is the grid width of the initial uniform configuration of the particles. For the 1D simulations, the initial conditions were sufficiently smooth and the approximation sufficiently high in resolution that the effect of regularization was marginal and not required to ensure stability. However, for the 2D simulation in Experiment 3, computational cost is a consideration and the simulation was run at too coarse a resolution to remain stable unless  $\hat{\alpha}$ was increased to quadruple the reference grid size. 

\subsection{Experiment 1: conservative properties of VFL}
Unless otherwise stated, we use 128 cells, where each cell is initially of the same size. The parameters of the simulations are provided in Table \ref{tab:exp_vfl} below. 

\vspace{1.0cm}

\begin{table}[!ht]
\begin{center}
\begin{tabular}{llclr}
\hline 
\vline & Parameter &\vline& Value&\vline\\
\hline
\vline& Number of cells $N$&\vline& 128&\vline\\
\vline &Time step &\vline& 0.01 &\vline\\
\vline& Domain length &\vline & 2$\pi$ & \vline\\
\vline& Initial conditions &\vline& $m_{\alpha}(t_0)=1+0.1 exp(-800(X_{\alpha}(t_0)-L/2)^2/ L^2)$
&\vline\\
\vline&&\vline& $U_{\alpha}(t_0)=0$ &\vline\\
\vline&&\vline& $V_{\alpha}(t_0)=0$ &\vline\\
\vline& $f_0$ & \vline & $2\pi$ & \vline\\
\vline& $g$ & \vline & $4\pi^2$ & \vline\\
\vline& $\hat{\alpha}$ & \vline & $dx$ & \vline\\
\hline
\end{tabular}
\end{center}
\caption{This Table lists the simulation parameters for experiment 1, which
investigates the conservative properties of the VFL method for 1D rotating shallow water.}
\label{tab:exp_vfl}
\end{table}

\begin{figure}[!ht] 
\begin{center}
\includegraphics[angle=0.0,width=0.8\textwidth, height=0.3\textheight]{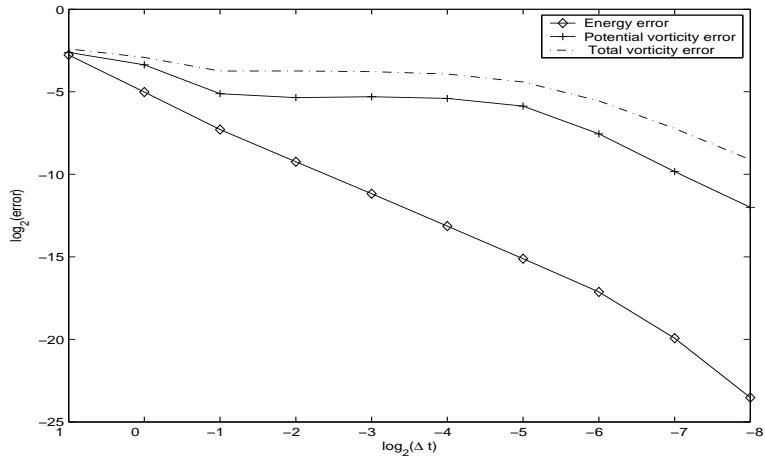}\\
\end{center}
\caption{This Figure shows the scaling laws of the mean energy error, potential vorticity and total vorticity with the size of the time step over $10^6$ time steps. The energy error scales to an order $\Delta t^2$ which is consistent with the order of the St\"ormer Verlet integrator. Potential vorticity (PV) and total vorticity (TV) error do not scale with the order of the integrator suggesting that the order of the error of the PV and TV are not governed by the integrator but by the VFL method itself. The simulation parameters are given in Table \ref{tab:exp_vfl}.}
\label{fig:error_scaling}
\end{figure}
\begin{figure}[!ht] 
\begin{center}
\includegraphics[angle=0.0,width=0.9\textwidth, height=0.3\textheight]{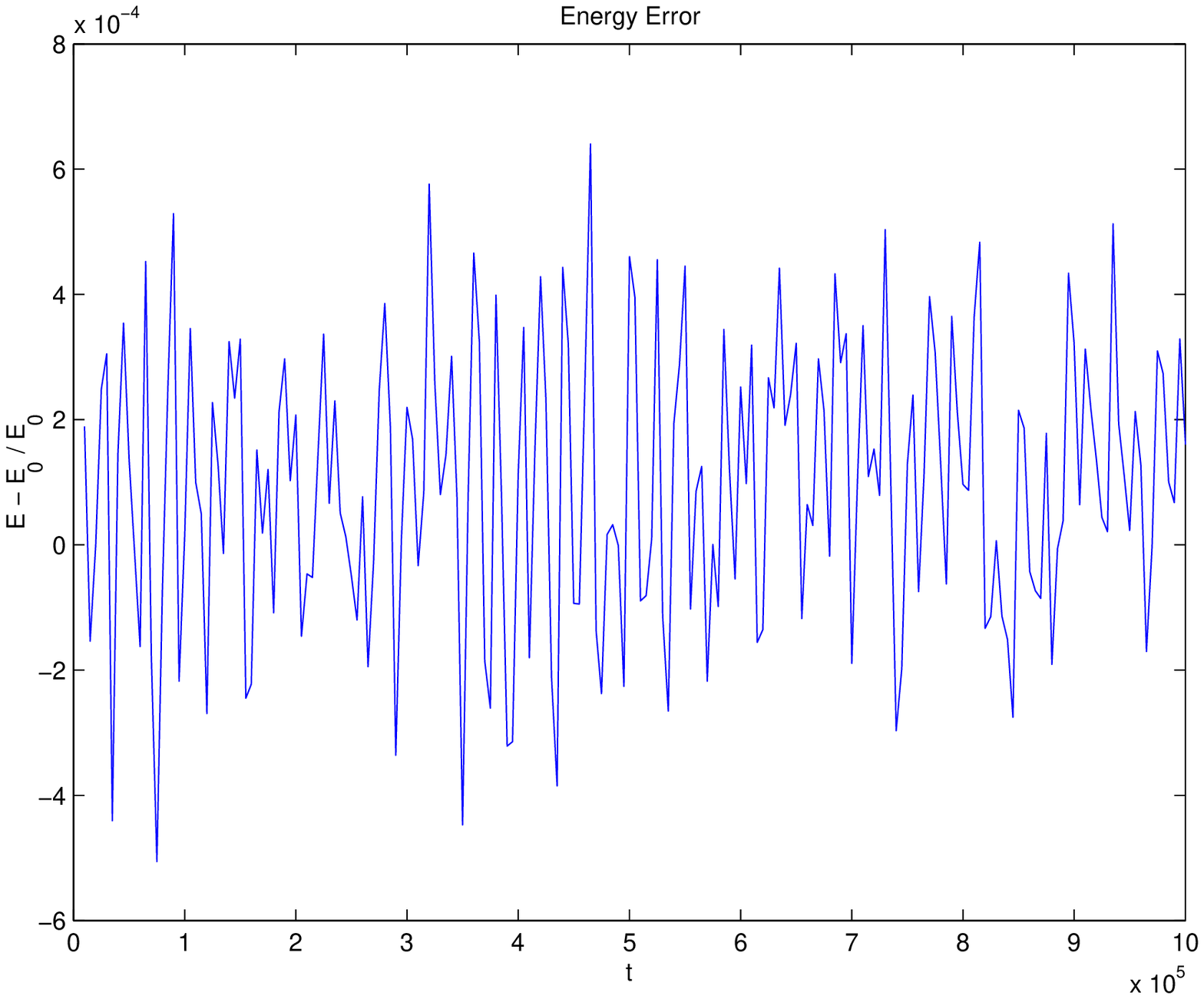}\\
\includegraphics[angle=0.0,width=0.9\textwidth,height=0.3\textheight]{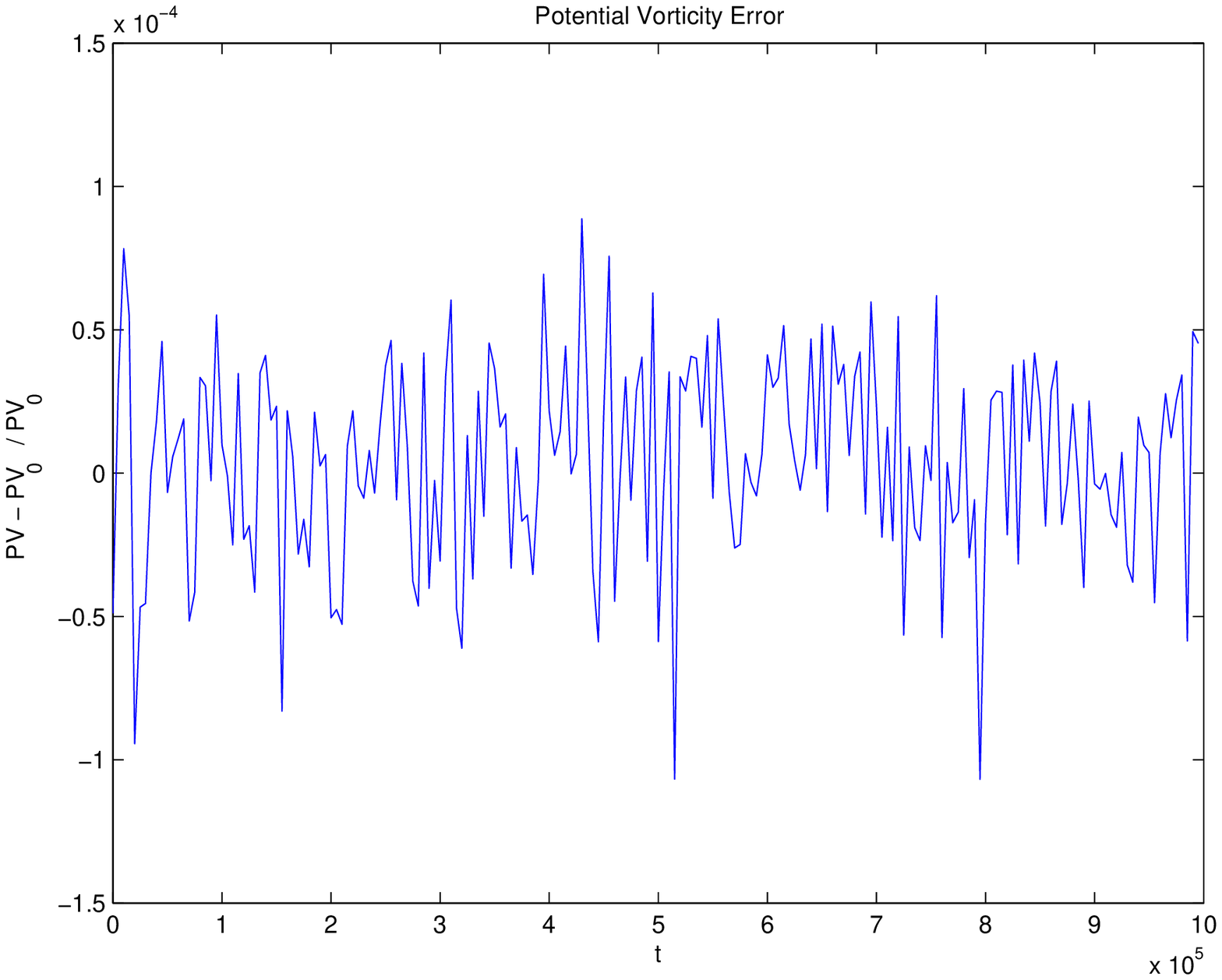}\\
\includegraphics[angle=0.0,width=0.9\textwidth,height=0.3\textheight]{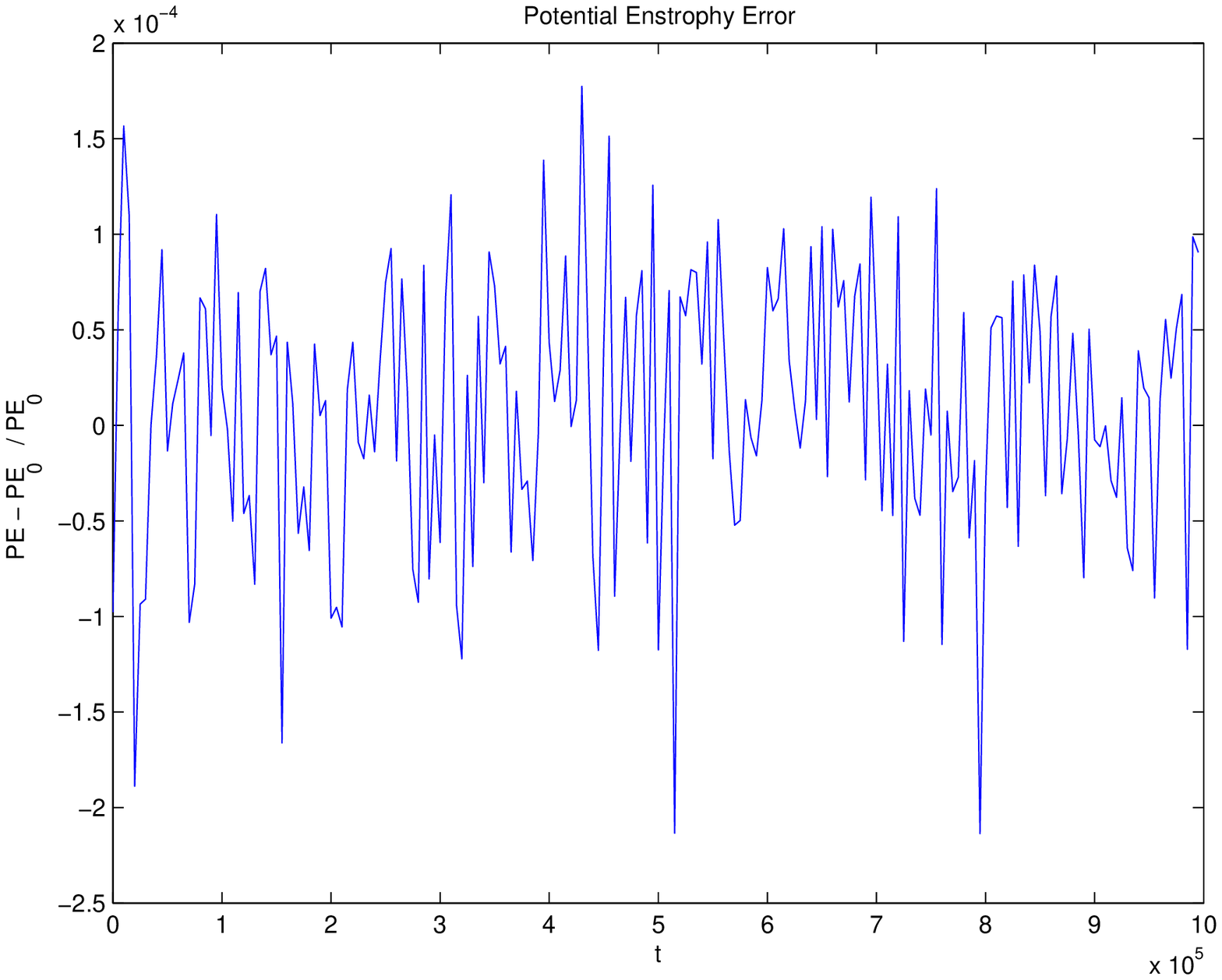}\\
\end{center}
\caption{These graphs show (from top to bottom) the relative energy, potential vorticity and potential enstrophy error over $10^6$ time steps. Each error does not drift which suggests that the VFL method is long-time stable and conservative. The simulation parameters are given in Table \ref{tab:exp_vfl}.}
\label{fig:vfl_time_err}
\end{figure}

\paragraph{\emph{Results}}  The following Figures show the relative energy, total potential vorticity (PV) and total potential enstrophy PE (sum of the square of the cellwise potential vorticity) errors. Figure \ref{fig:error_scaling} shows the scaling of the energy, PV and total vorticity error with the size of the time step. All three quantities montonically decrease as the time step decreases. The error in the energy scales with the order of integrator, $\mathcal{O}(\Delta t^2)$. where as the potential and total vorticity do not scale with this order.  This suggests that the VFL method introduces error in the discrete vorticity.  Figure \ref{fig:vfl_time_err} shows (from top to bottom) the relative energy, potential vorticity and potential enstrophy error over $10^6$ time steps, with a time step of $0.01$. This number of time steps is approximately equal to the number of time steps required to conduct a 100 year climate change simulation. Each error does not drift which suggests that the VFL method is long-time stable and conservative.

\newpage

\subsection{Experiment 2: Geostrophic adjustment}

This experiment uses the VFL method to simulate the mechanism of geostrophic adjustment of shallow water in a plane in which the layer depth and velocity are independent of the y-direction (meridional direction)

\begin{equation} \label{eqn:geostrophicBalance}
\dot{U}=-g \nabla_X h +f_0V,\quad
\dot{V}=-f_0U,\quad
\dot{X}=U,
\end{equation}

Given a geostrophically balanced layer depth $H^0$, meridional velocity $V^0=\frac{g}{f_0} \nabla_X h$ and horizontal velocity $U^0=0$, we perturb the layer depth $h=H^0
+ h'$ by a smooth, domain centred, gaussian function $h'$ with support $L'$ less than the Rossby deformation
radius $L_D$. The dynamics, governed by equations \eqref{eqn:geostrophicBalance}, are observed to exhibit gravity waves which propagate energy and momentum away from the source leaving behind geostrophically balanced flow.

We
verify in this numerical experiment that this mechanism only occurs if the scale of layer depth perturbation
is smaller than the Rossby deformation radius by performing two simulations,
one in which the scale of perturbation is larger and one is which it is smaller
than the Rossby deformation radius. The
parameters for this simulations are provided in Table \ref{tab:exp_vfl_2}.

\begin{table}[!ht]
\begin{center}
\begin{tabular}{llclr}
\hline 
\vline & Parameter &\vline& Value&\vline\\
\hline
\vline& Number of cells $N$&\vline& 512&\vline\\
\vline &Time step &\vline& 0.01 &\vline\\
\vline& Domain length &\vline & 2$\pi$ & \vline\\
\vline& Initial conditions &\vline& $m_{\alpha}(t_0)=1+\sum_{i=1}^2 0.01 exp(-\beta_i(X_{\alpha}(t_0)-L/2)^2/ L^2)$
&\vline\\
\vline&&\vline& $U_{\alpha}(t_0)=0$ &\vline\\
\vline&&\vline& $V_{\alpha}(t_0)=-\frac{g}{f_0}\text{grad}(\bar{H}^0_{\alpha})$ &\vline\\
\vline& $f_0$ & \vline & $2\pi$ & \vline\\
\vline& $g$ & \vline & $4\pi^2$ & \vline\\
\vline& $H_0$ & \vline & $1$ & \vline\\
\vline& $L_D$ & \vline & $2\pi$ & \vline\\
\vline& $\hat{\alpha}$ & \vline & $dx$ & \vline\\
\hline
\end{tabular}
\end{center}
\caption{This Table lists the simulation parameters for experiment 2 which
models geostrophic adjustment when $\beta_1=10 (L'<<L_D)$ and $\beta_2=1000 (
L'>>L_D)$. $\bar{H}^0_{\alpha}$ denotes the discrete steady state profile, from which the meridional velocity is initialized.}
\label{tab:exp_vfl_2}
\end{table}

\begin{figure}[!ht] 
\begin{center}
\includegraphics[angle=0.0,width=0.7\textwidth]{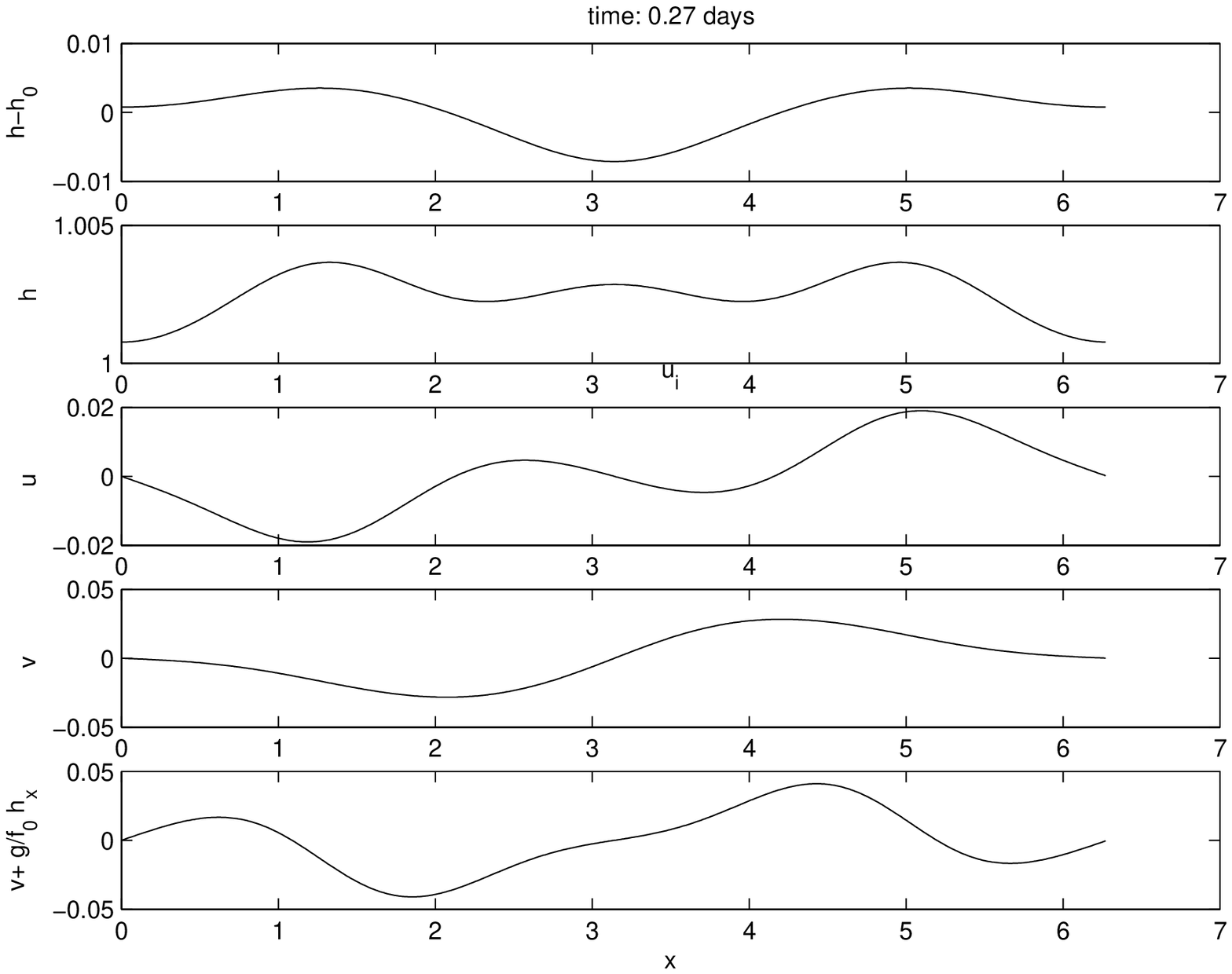}
\end{center}
\caption{These graphs show a snapshot of the shallow water taken at time
$t=0.27$ "days" (1 "day" is a rotational unit) for which the initial scale
of perturbation $L'>>L_D$, where $L_D$ is the Rossby deformation radius. The top two graphs show the layer
depth perturbation and layer depth respectively. The third and fourth graphs
show the horizontal and meridional velocities.  The bottom graph shows the
difference of the meridional velocity with the layer depth gradient. Note that the source region is
\emph{not} restored to a geostrophically balanced state because geostrophic adjustment does not occur. The simulation parameters for this experiment are given in Table \ref{tab:exp_vfl_2} in which $\beta=1000$.}
\label{fig:vfl_rossby_1}
\end{figure}

\begin{figure}[!ht] 
\begin{center}
\includegraphics[angle=0.0,width=0.6\textwidth]{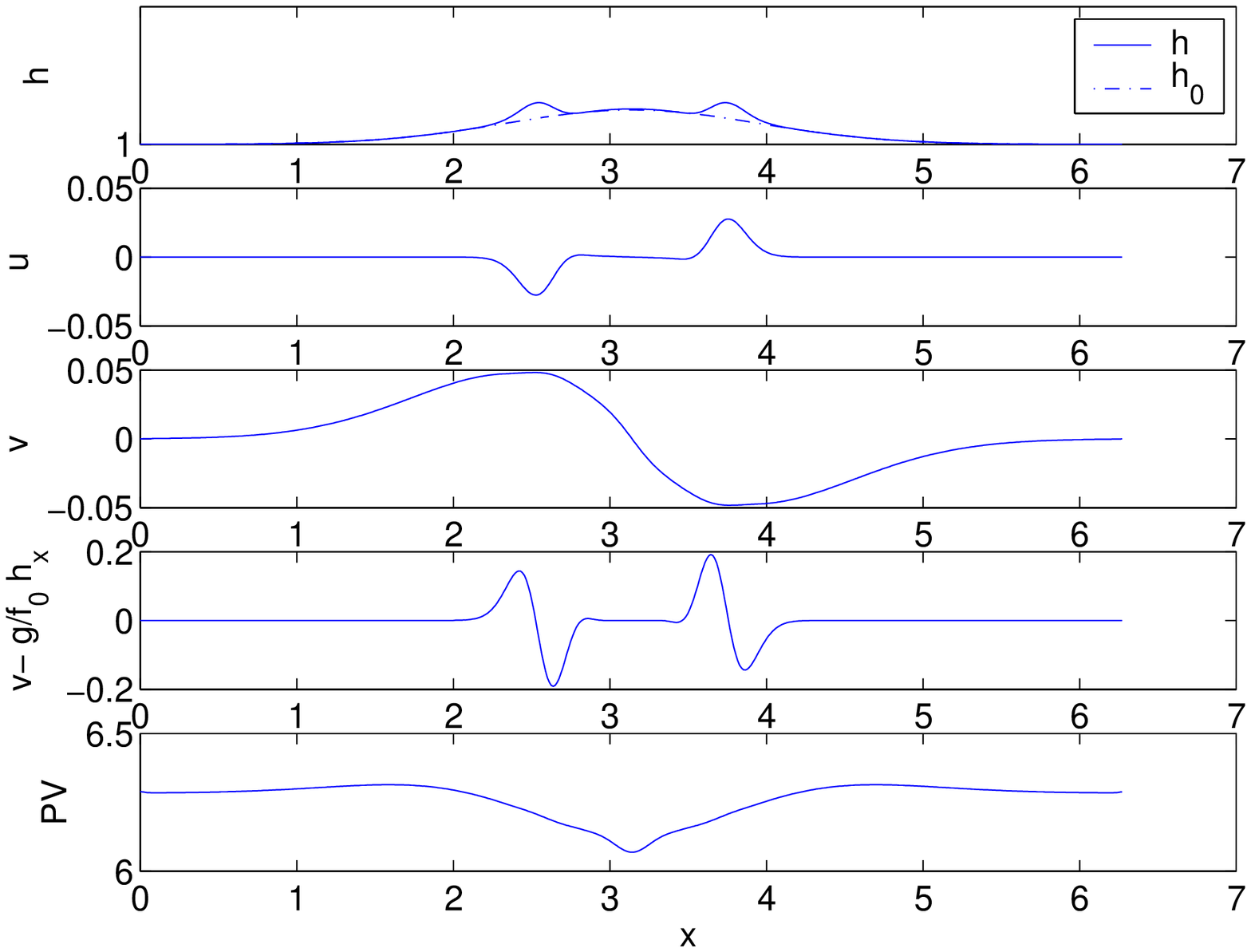}\\
\includegraphics[angle=0.0,width=0.6\textwidth]{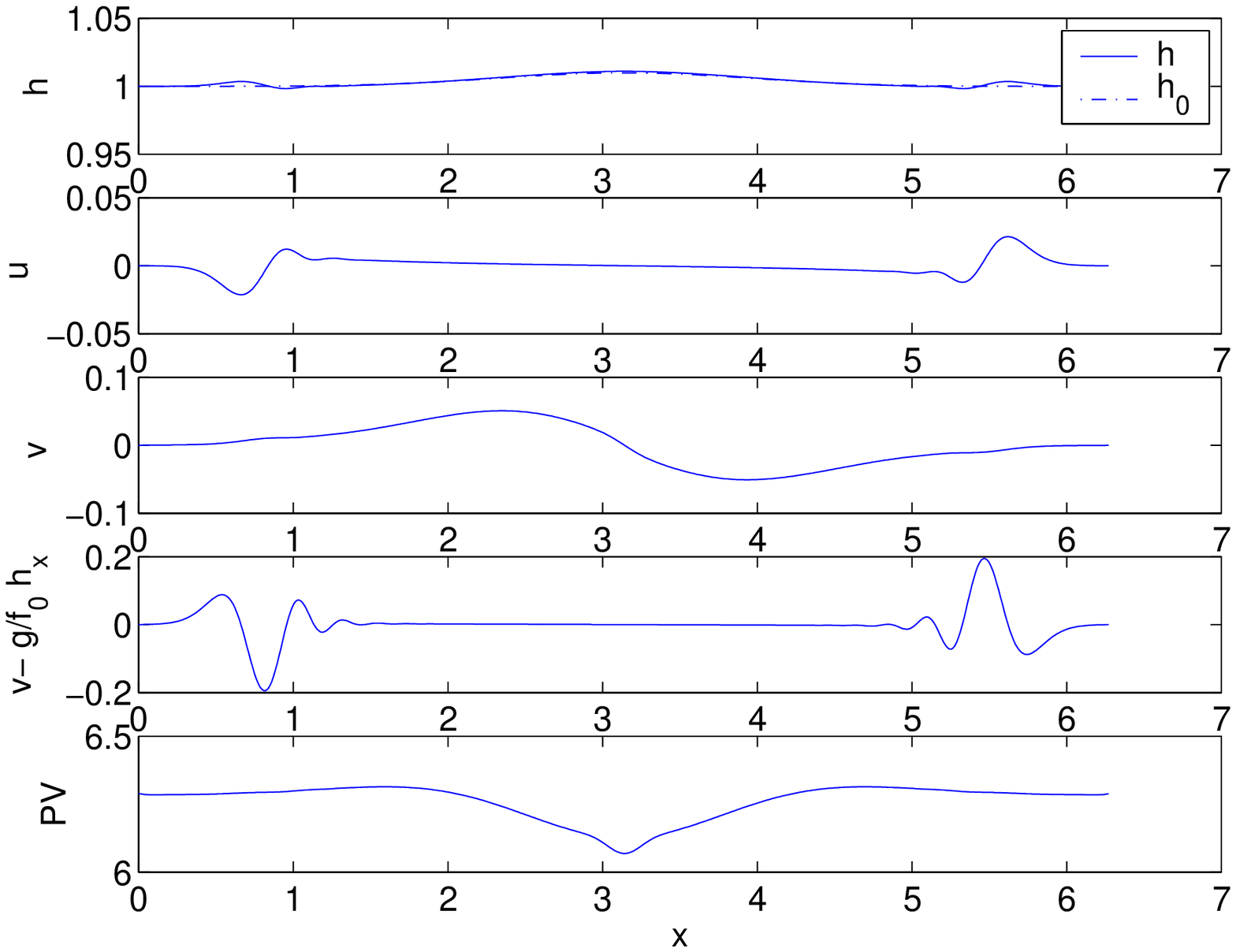}\\
\end{center}
\caption{These graphs show snapshots of the shallow water taken at times
$t=\{0.2,0.5\}$ "days" (1 "day" is a rotational unit) for which the initial scale
of perturbation $L'<<L_D$, where $L_D$ is the Rossby deformation radius. The top graph show the layer
depth perturbation and layer depth respectively. The third and fourth graphs
show the horizontal and meridional velocities.  The bottom graph shows
the difference of the meridional
velocity with the layer depth gradient. This difference is zero when the flow is geostrophically balanced. Note that the source region does
restore to a geostrophically balanced state because gravity waves propagate
away energy and momentum to restore the balance. The simulation parameters for this experiment are given in Table \ref{tab:exp_vfl}, in which $\beta=50$.}
\label{fig:vfl_rossby_2}
\end{figure}

\paragraph{\emph{Results}} Figures \ref{fig:vfl_rossby_1} and \ref{fig:vfl_rossby_2} compare two different rotating shallow water regimes distinguished by the scale of perturbation of the layer depth. In the first case, the layer depth is perturbed on the scale of the Rossby radius $L_R=1$ and in the second, on a scale smaller than the Rossby radius. In the latter case, the sequence of layer depth graphs (shown at increasing simulation times) in Figure \ref{fig:vfl_rossby_2}
show the gravity waves which propagate from the source. The bottom profile of each of these graphs in this Figure also shows that a geostrophically balanced region is recovered in the region of the source.

\begin{figure}[!ht] 
\begin{center}
\includegraphics[angle=0.0,width=0.7\textwidth]{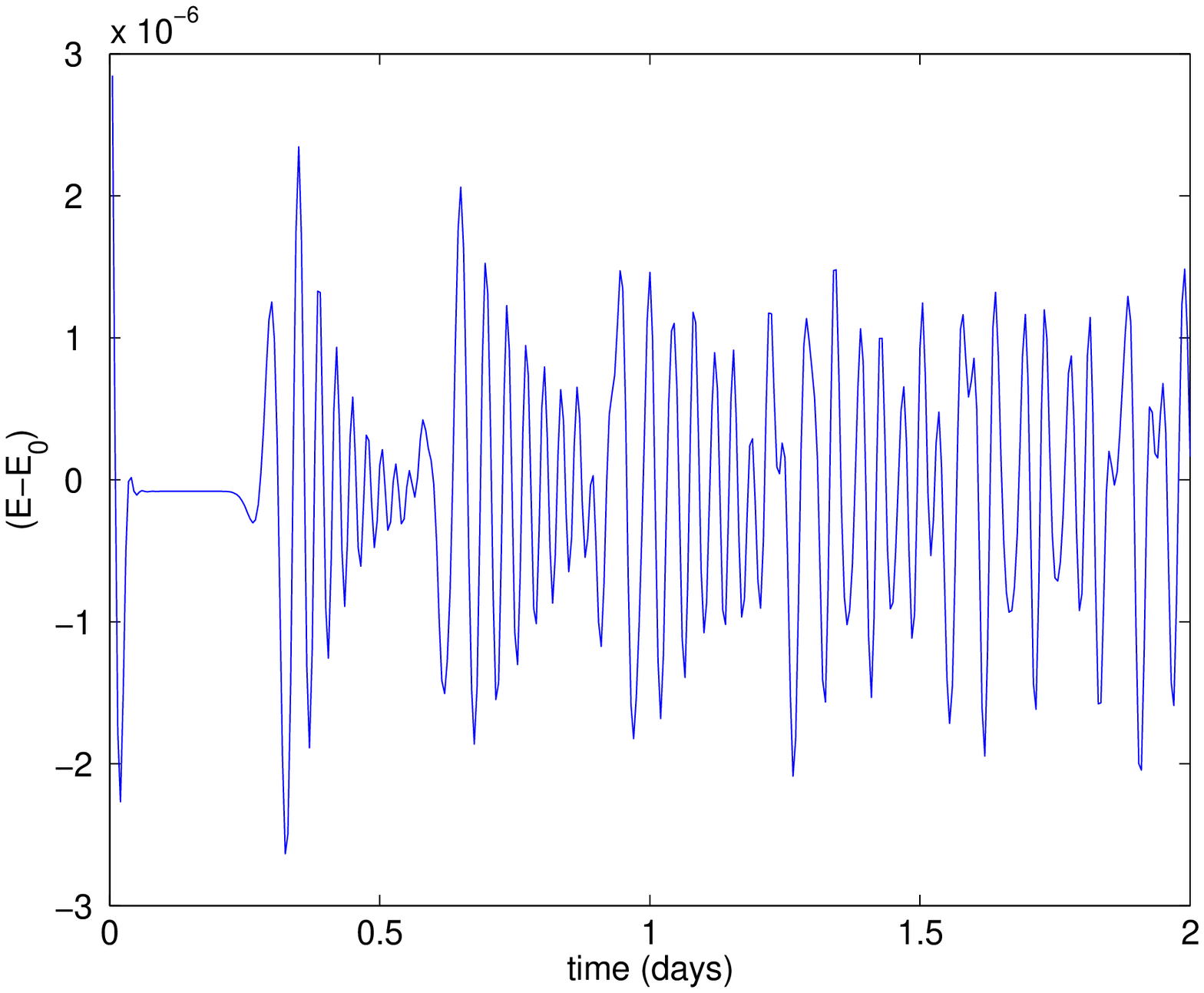}\\
\includegraphics[angle=0.0,width=0.7\textwidth]{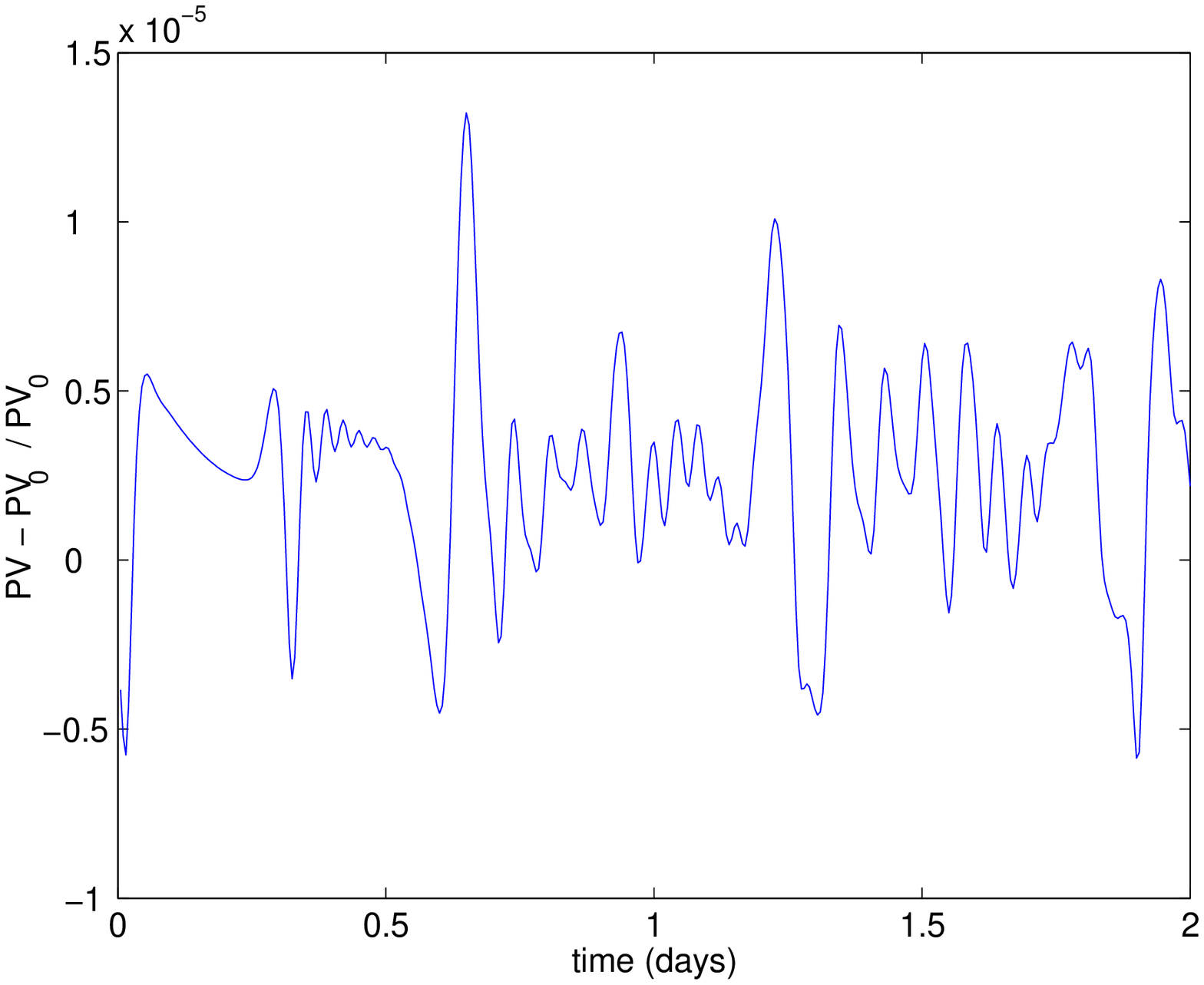}\\
\end{center}
\caption{These graphs show the relative energy (top) and potential vorticity (bottom) error profiles of shallow water over 2 "days" (1 "day" is a rotational unit) for which the initial scale
of perturbation is $L'<<L_D$, where $L_D$ is the Rossby deformation radius.
Both profiles exhibit high frequency oscillations arising from gravity waves but do not exhibit drift. There are also jumps in the profiles with an approximate period of $0.3$ days.
These jumps occur when the gravity waves collide (recall that the domain
is periodic).  The simulation parameters for this experiment are given in Table \ref{tab:exp_vfl_2} in which $\beta=50$.}
\label{fig:vfl_rossby_4}
\end{figure}

\begin{figure}
\begin{center}
\includegraphics[angle=0.0,width=0.6\textwidth]{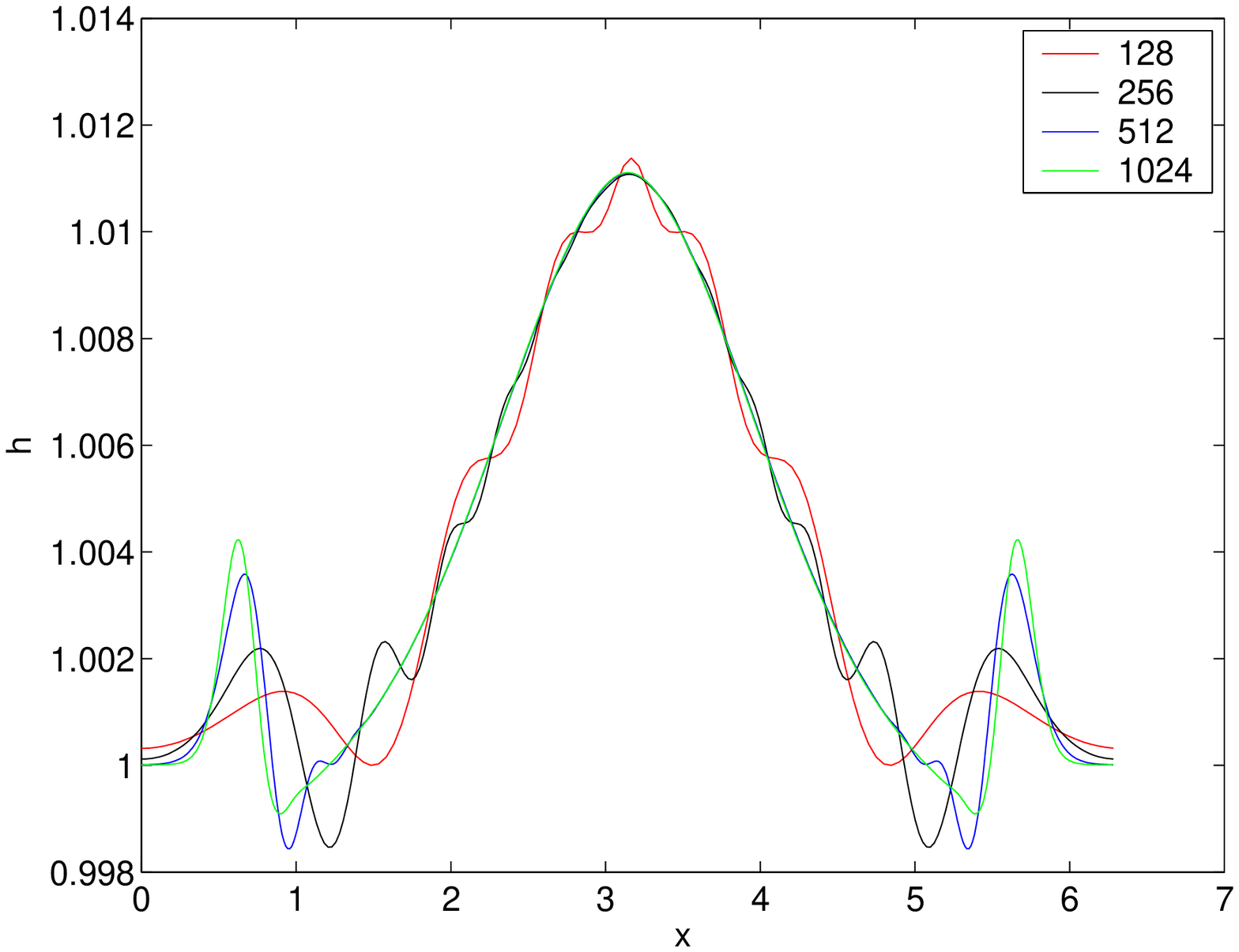}\\
\end{center}
\caption{This graph shows the layer depth at $0.4$ days for various number of grid points. In the central region of the domain, the layer depth has reached a steady state (geostrophic balance). The graph shows how the discrete approximation of the layer depth converges to the steady state. The outer regions of the domain exhibit gravity waves which have propagated away from the source. This outer region is not in geostrophic balance after 0.4 days.}
\label{fig:chpt_convergence}
\end{figure}

Figure \ref{fig:vfl_rossby_4} show the corresponding
energy and potential vorticity errors over 2 rotation units (i.e. 2 days).
We firstly observe that the profiles exhibit high frequency oscillations arising from gravity waves but do not exhibit drift.
There are also jumps in the profiles with an approximate period of $0.3$ days.
These jumps occur when the gravity waves collide (recall that the domain
is periodic). Figure \ref{fig:chpt_convergence} shows how the discrete approximation of the layer depth, over the central region of the domain, converges to the steady state after
 $0.4$ days for various number of grid points. The outer regions of the domain exhibit gravity waves which have propagated away from the source. 

\clearpage
\subsection{Experiment 3: Vortex pair in 2D rotating shallow water flow}
This experiment shows the motion of a vortex
pair in a purely rotational shallow water flow over a doubly periodic domain without bottom topography. The flow field is initially smooth and geostrophic balanced - the layer thickness and horizontal velocity are radially symmetric about the center of the domain and the meridional velocity is anti-symmetric about the lines $x=L/2$ and $y=L/2$. The simulation parameters are given in Table \ref{tab:2d_conditions}. 

\begin{figure} 
\begin{center}
\includegraphics[angle=0.0,width=0.3\textwidth]{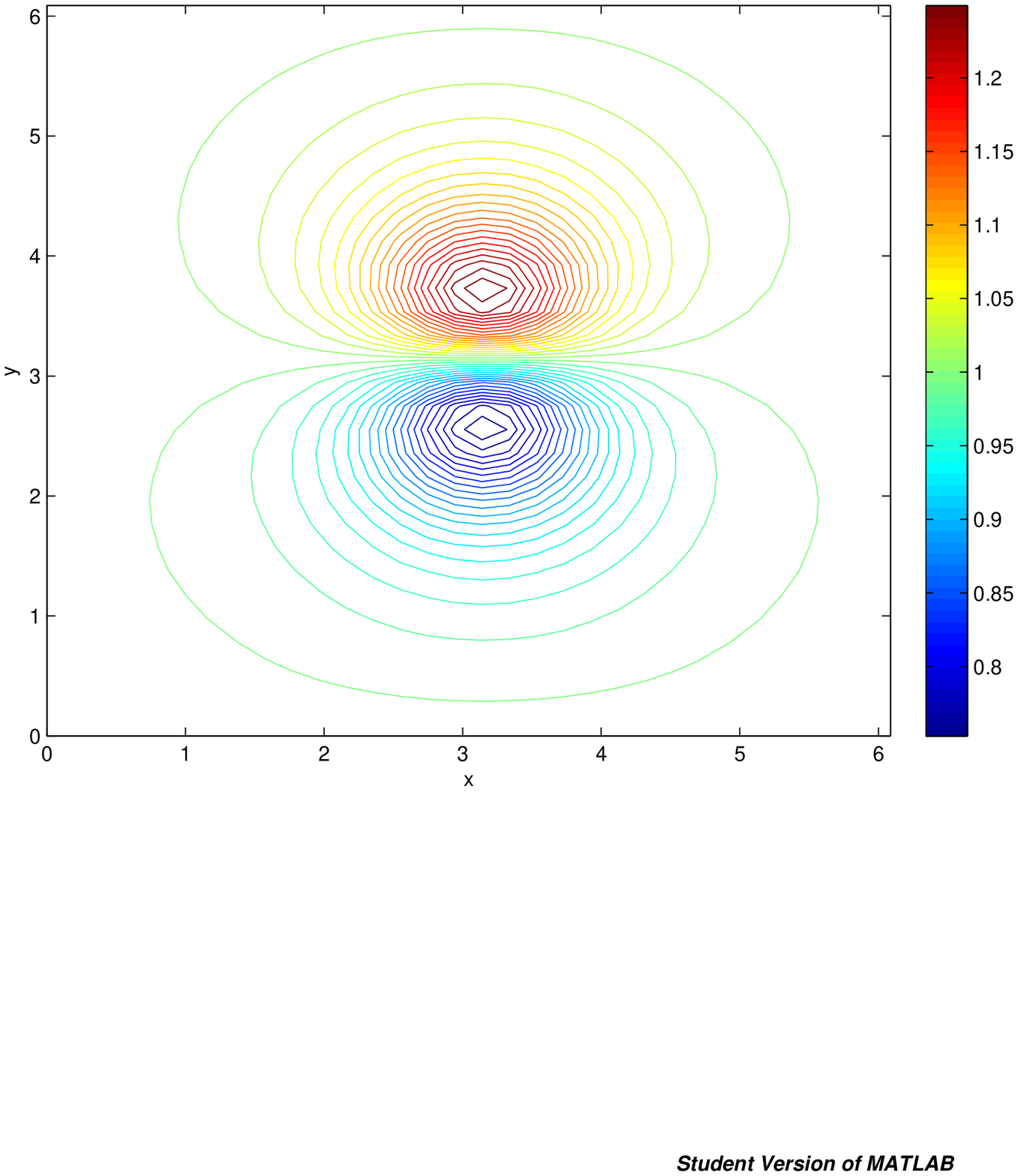}
\includegraphics[angle=0.0,width=0.3\textwidth]{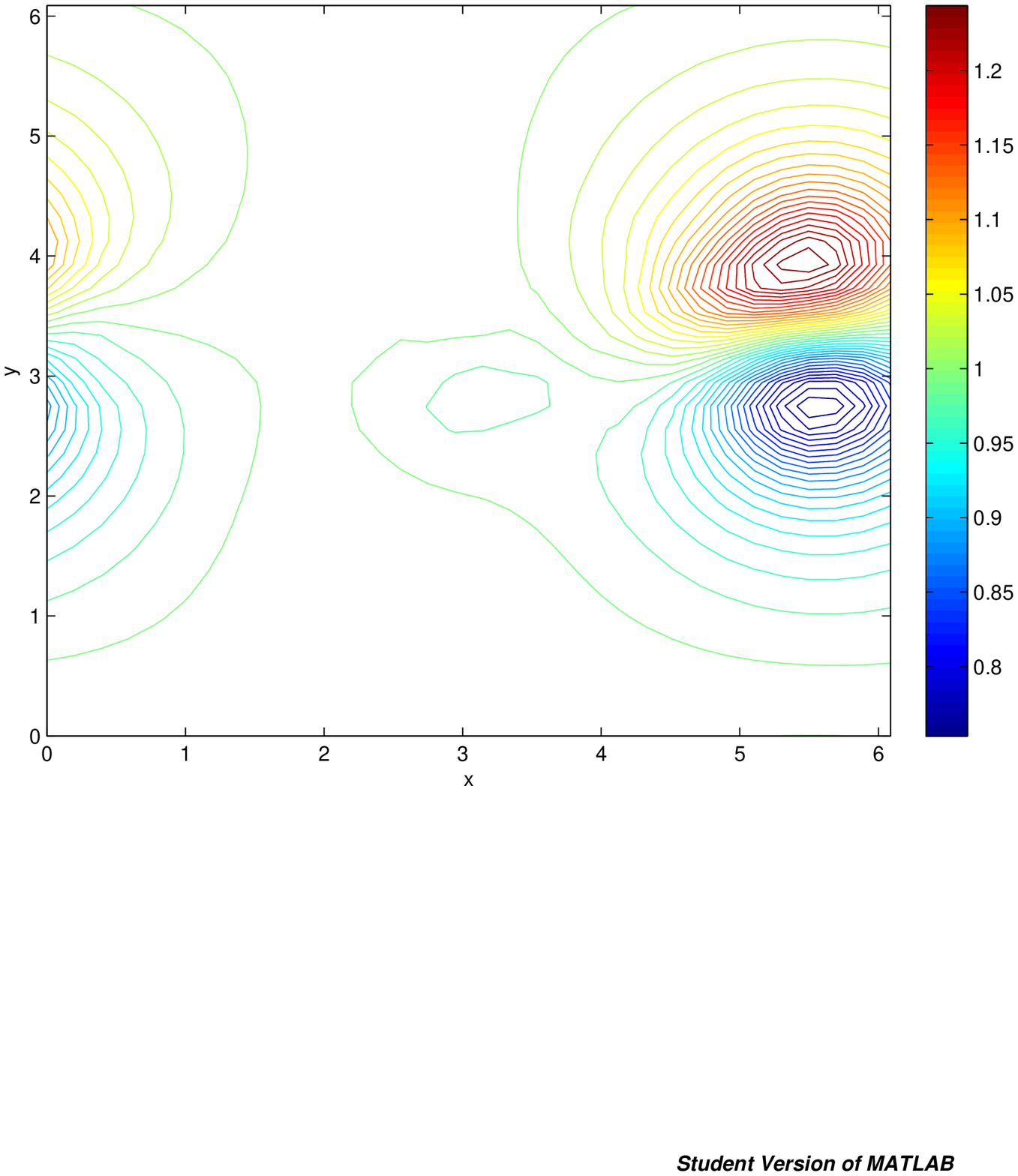}\\
\includegraphics[angle=0.0,width=0.3\textwidth]{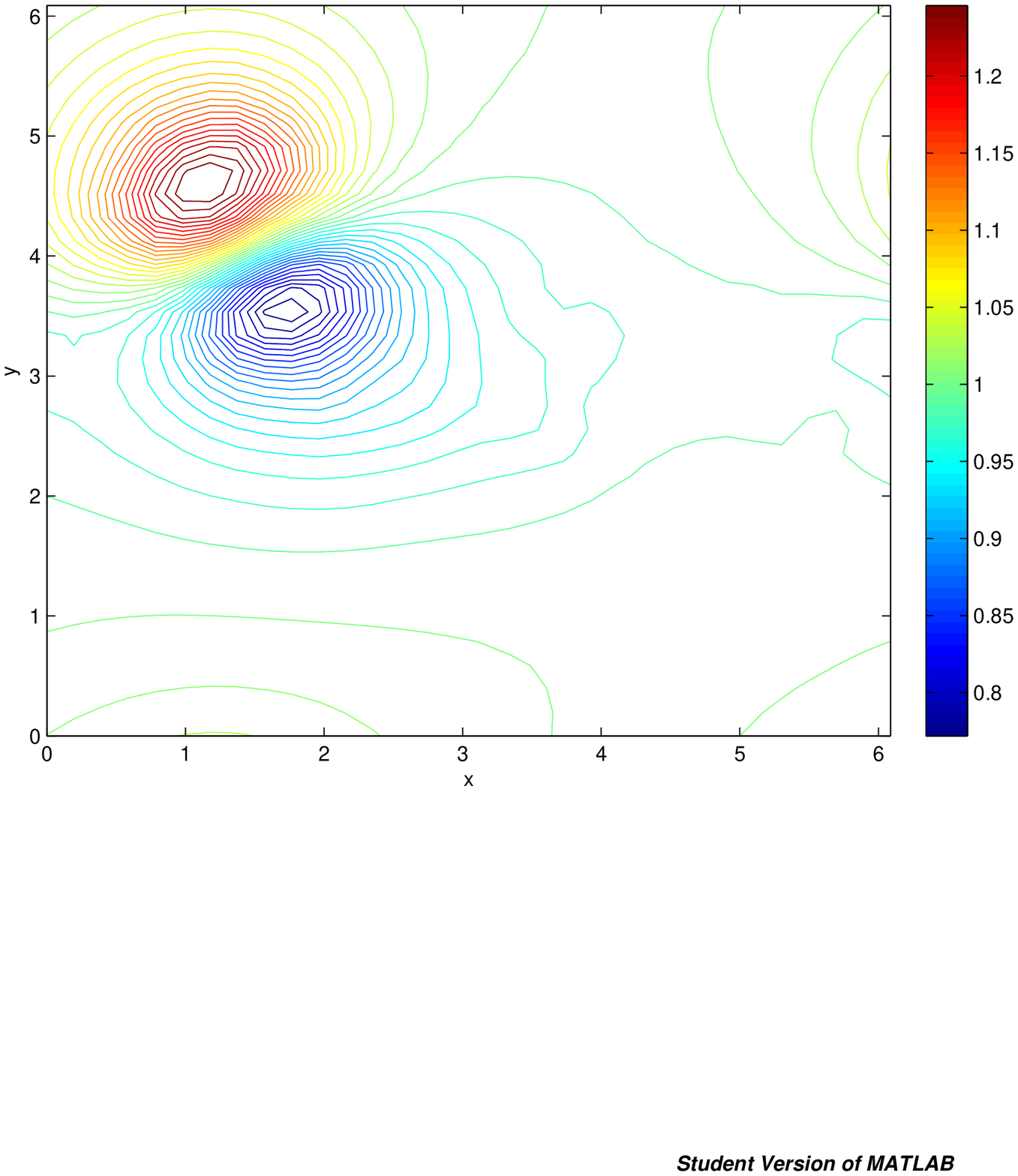}
\includegraphics[angle=0.0,width=0.3\textwidth]{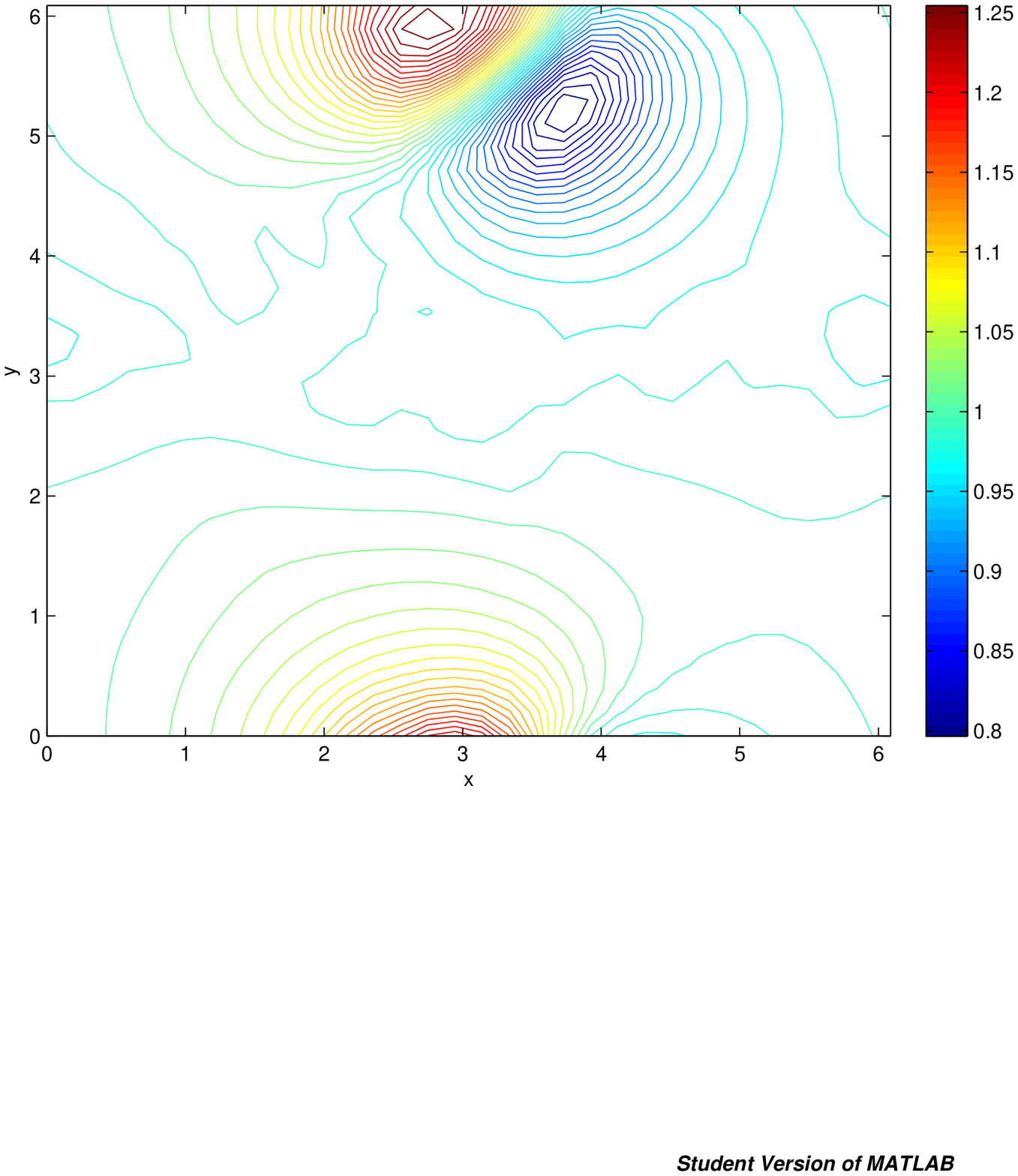}\\
\includegraphics[angle=0.0,width=0.3\textwidth]{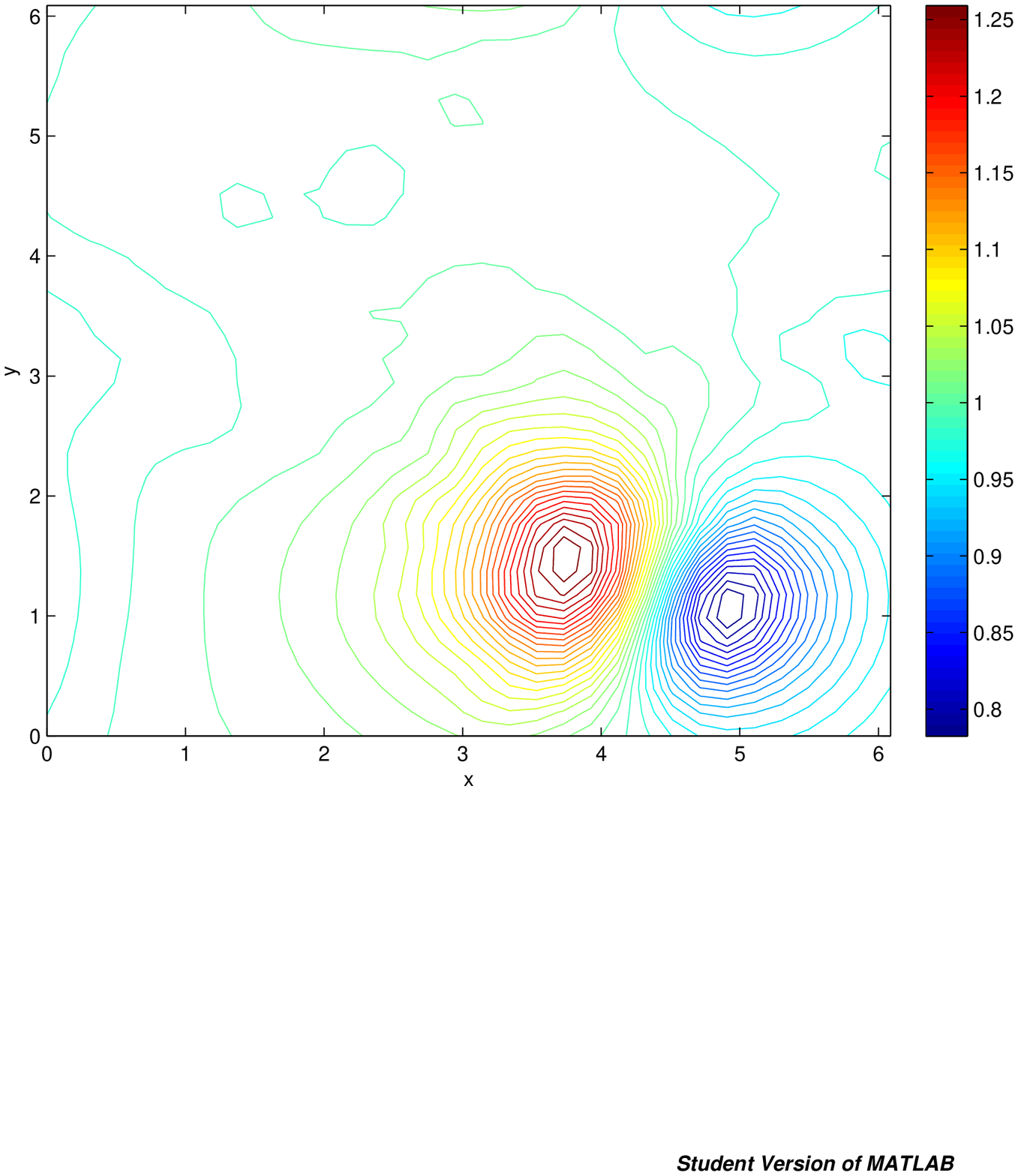}
\includegraphics[angle=0.0,width=0.3\textwidth]{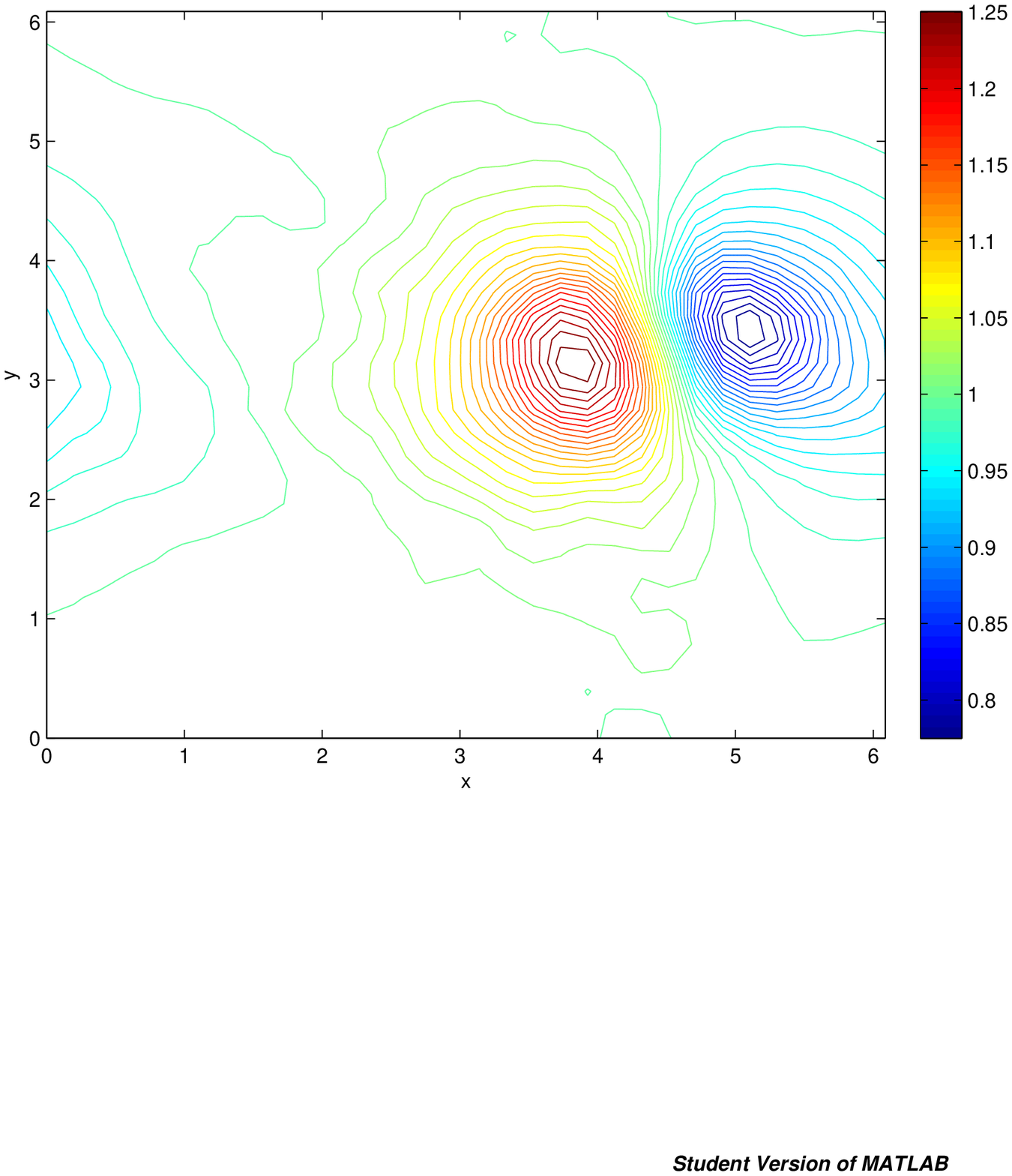}\\
\includegraphics[angle=0.0,width=0.3\textwidth]{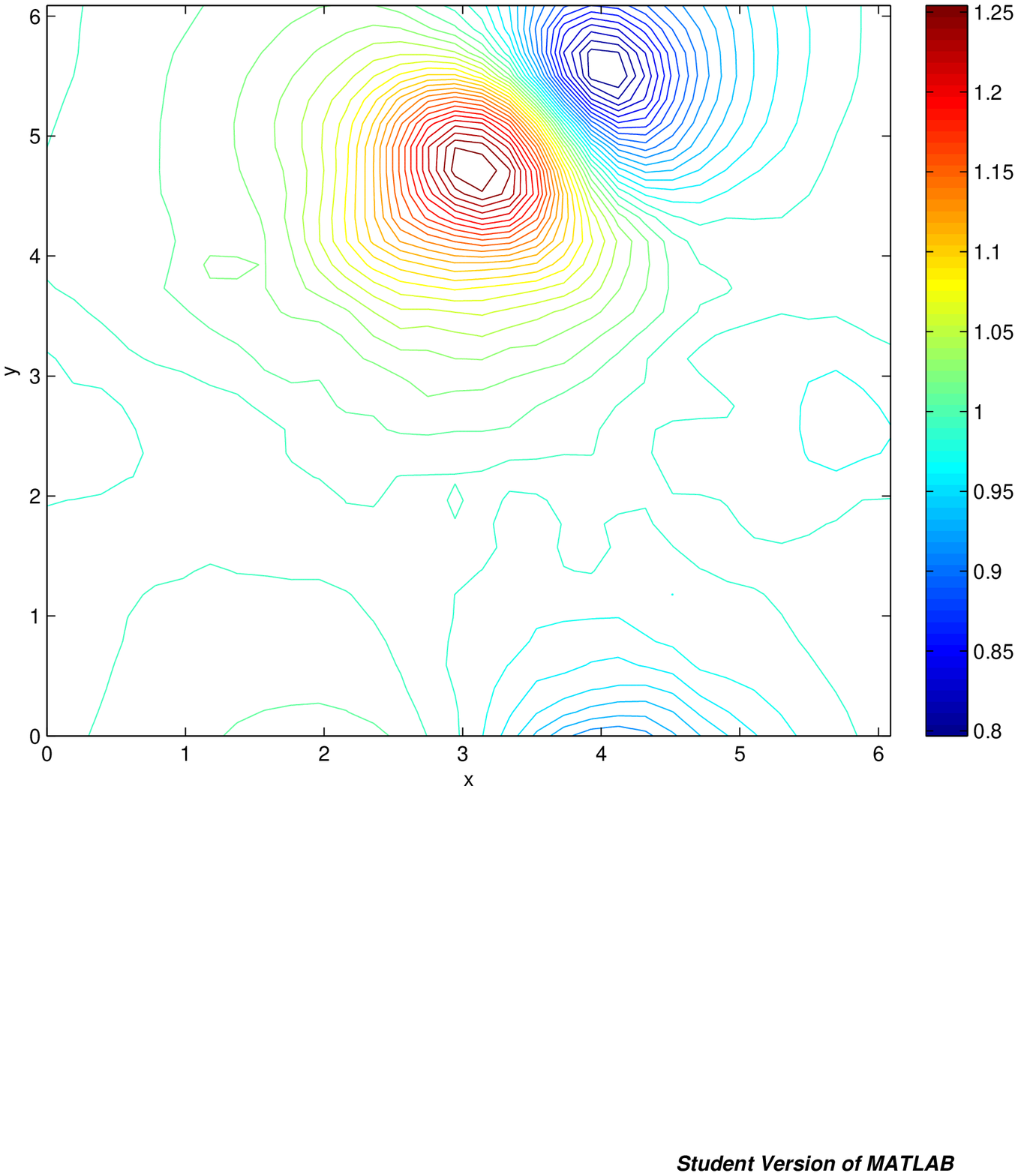}
\includegraphics[angle=0.0,width=0.3\textwidth]{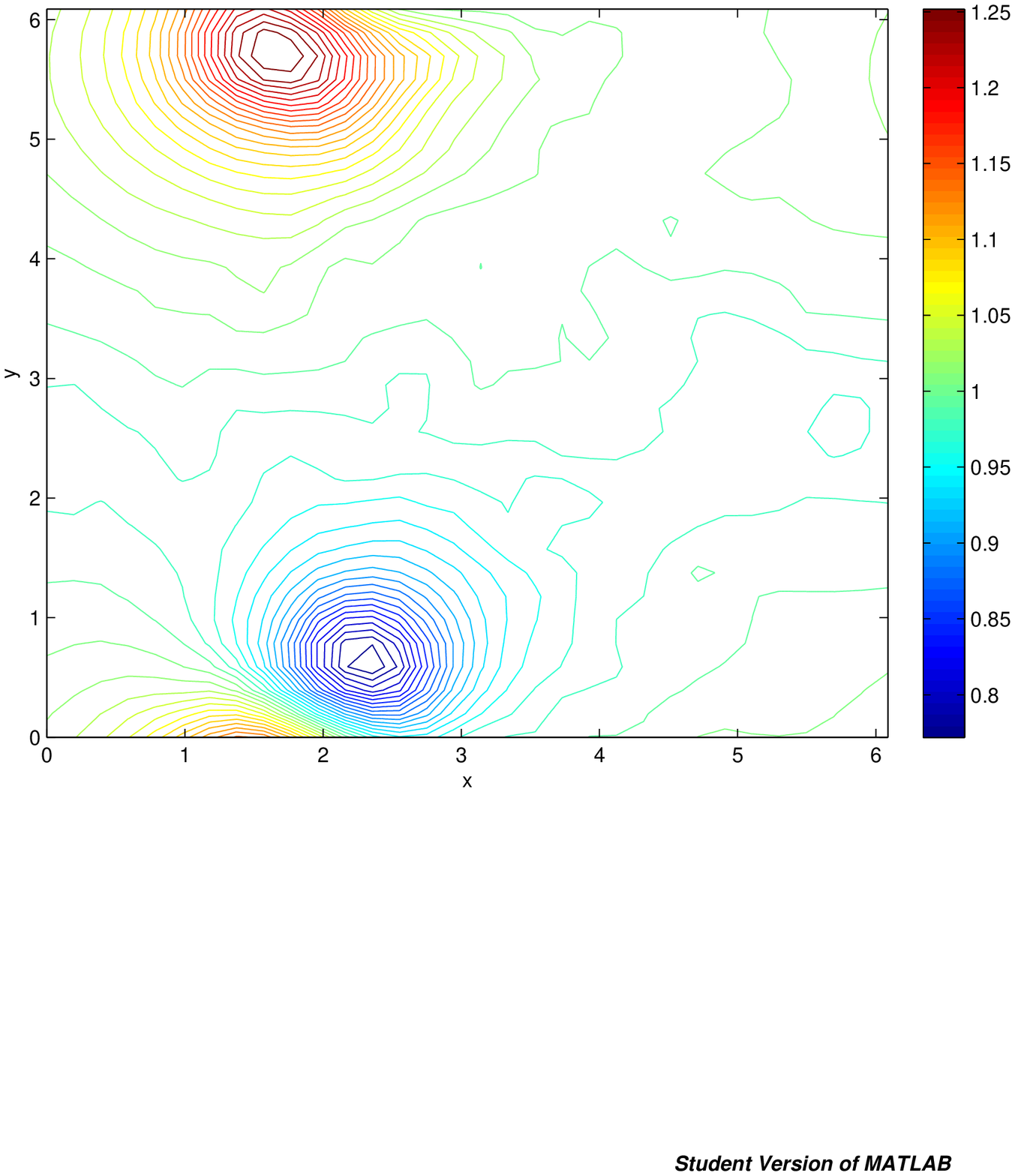}\\
\end{center}
\caption{This graph shows the layer depth at increments of $40$ days during simulation of a vortex pair in rotating shallow water over a doubly periodic domain. The simulation parameters for this experiment are given in Table \ref{tab:2d_conditions}.}
\label{fig:vfl_vortex}
\end{figure}

\begin{figure} 
\begin{center}
\includegraphics[angle=0.0,width=0.3\textwidth]{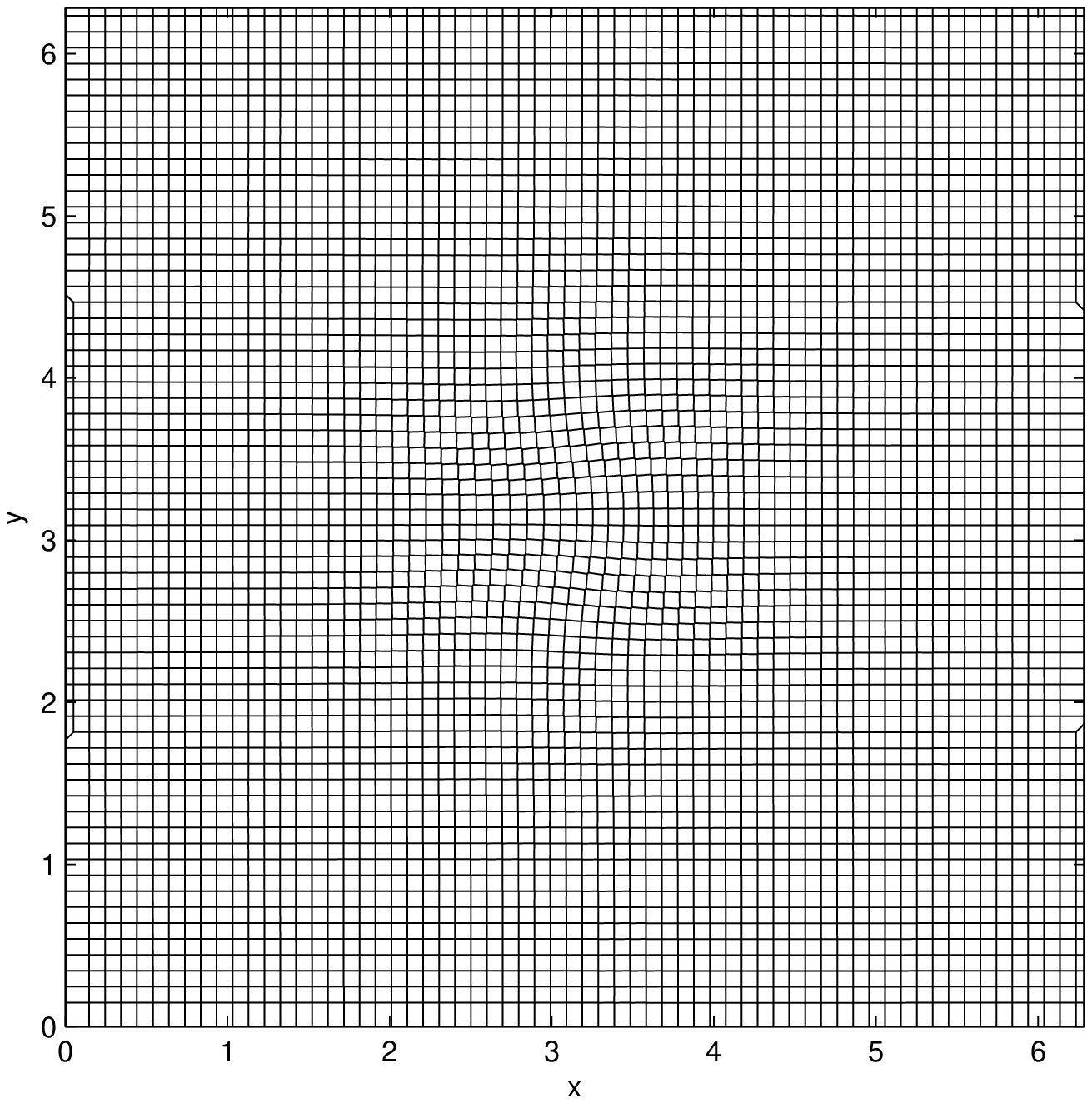}
\includegraphics[angle=0.0,width=0.3\textwidth]{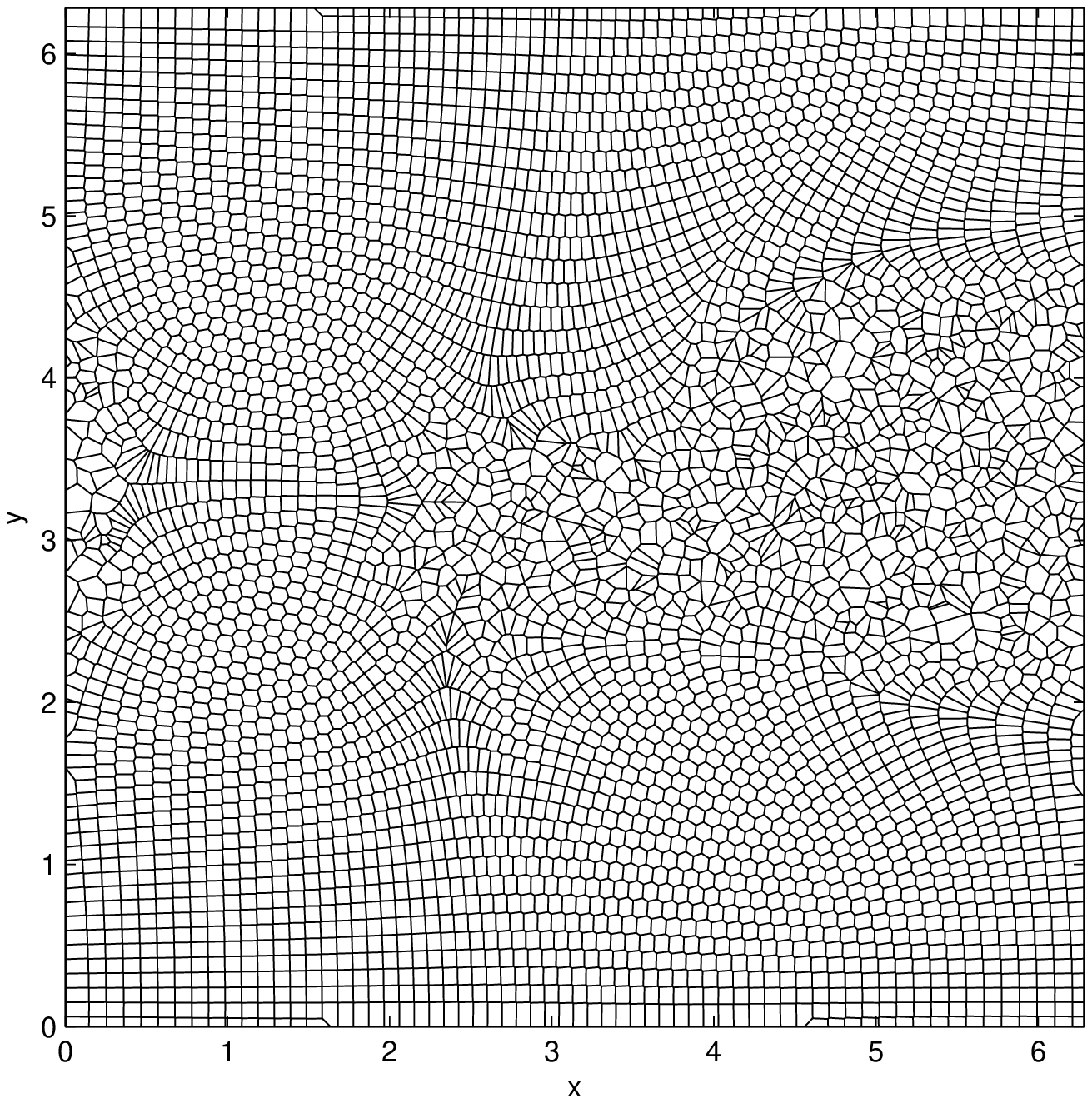}\\
\includegraphics[angle=0.0,width=0.3\textwidth]{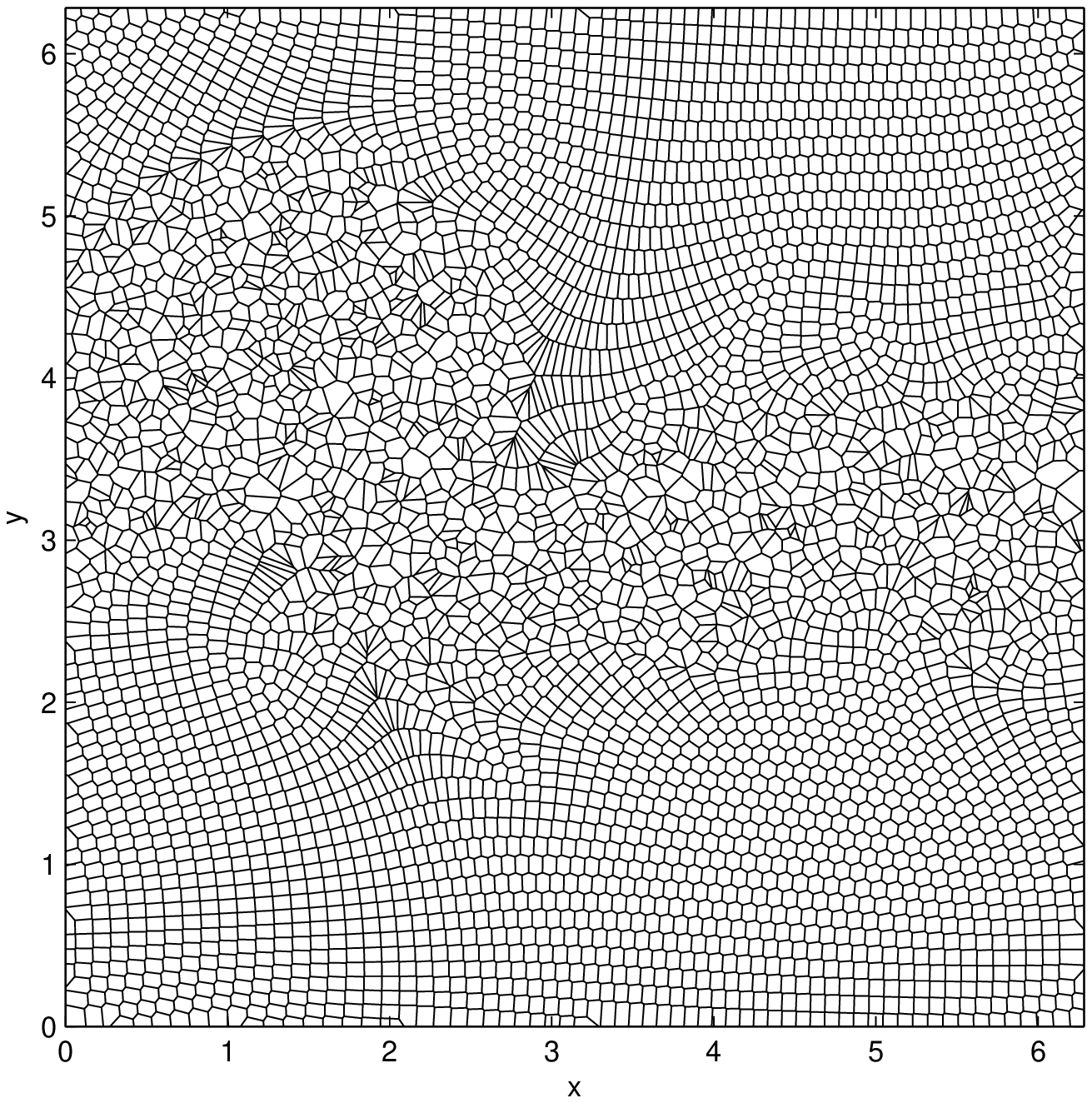}
\includegraphics[angle=0.0,width=0.3\textwidth]{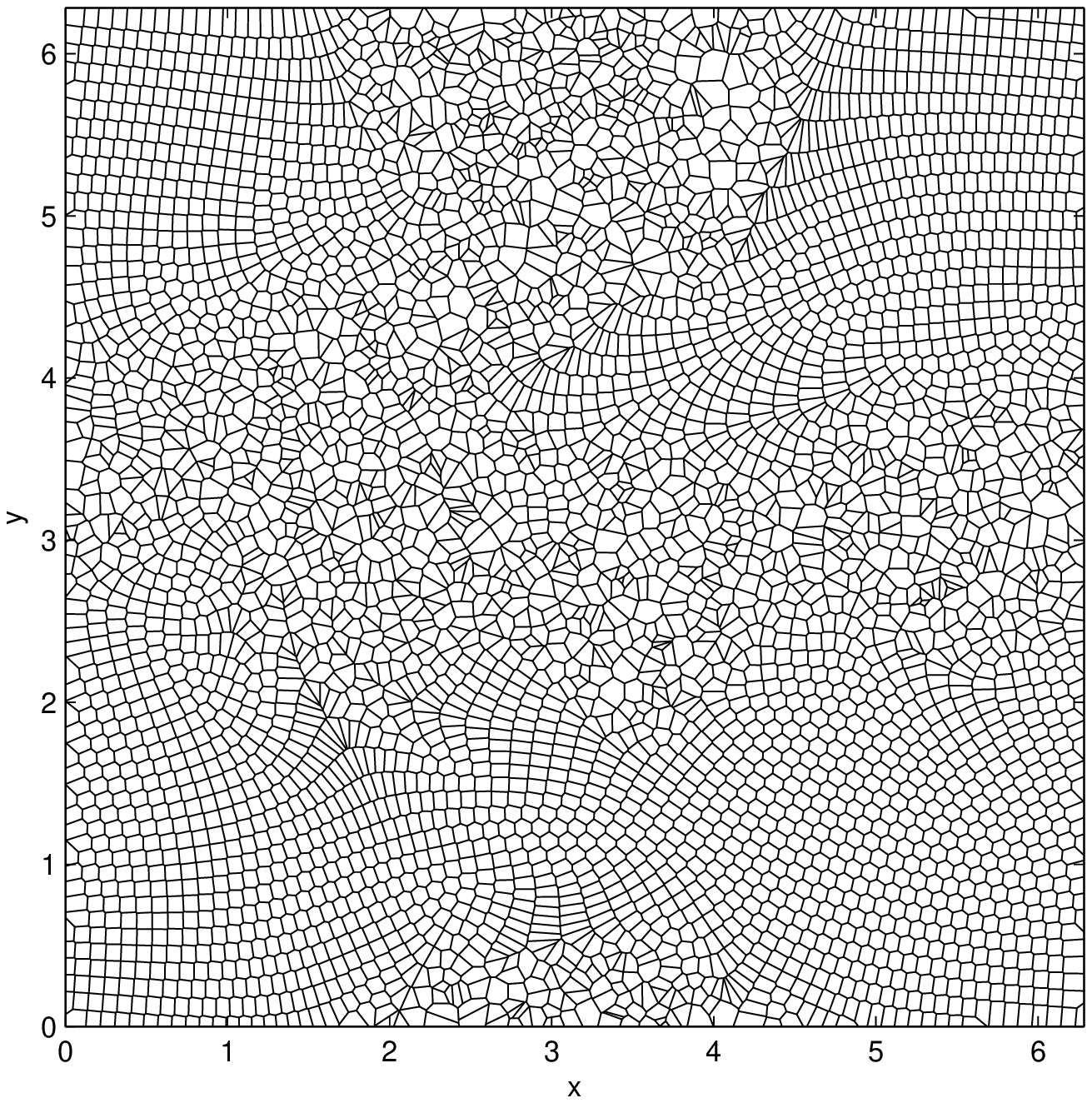}\\
\includegraphics[angle=0.0,width=0.3\textwidth]{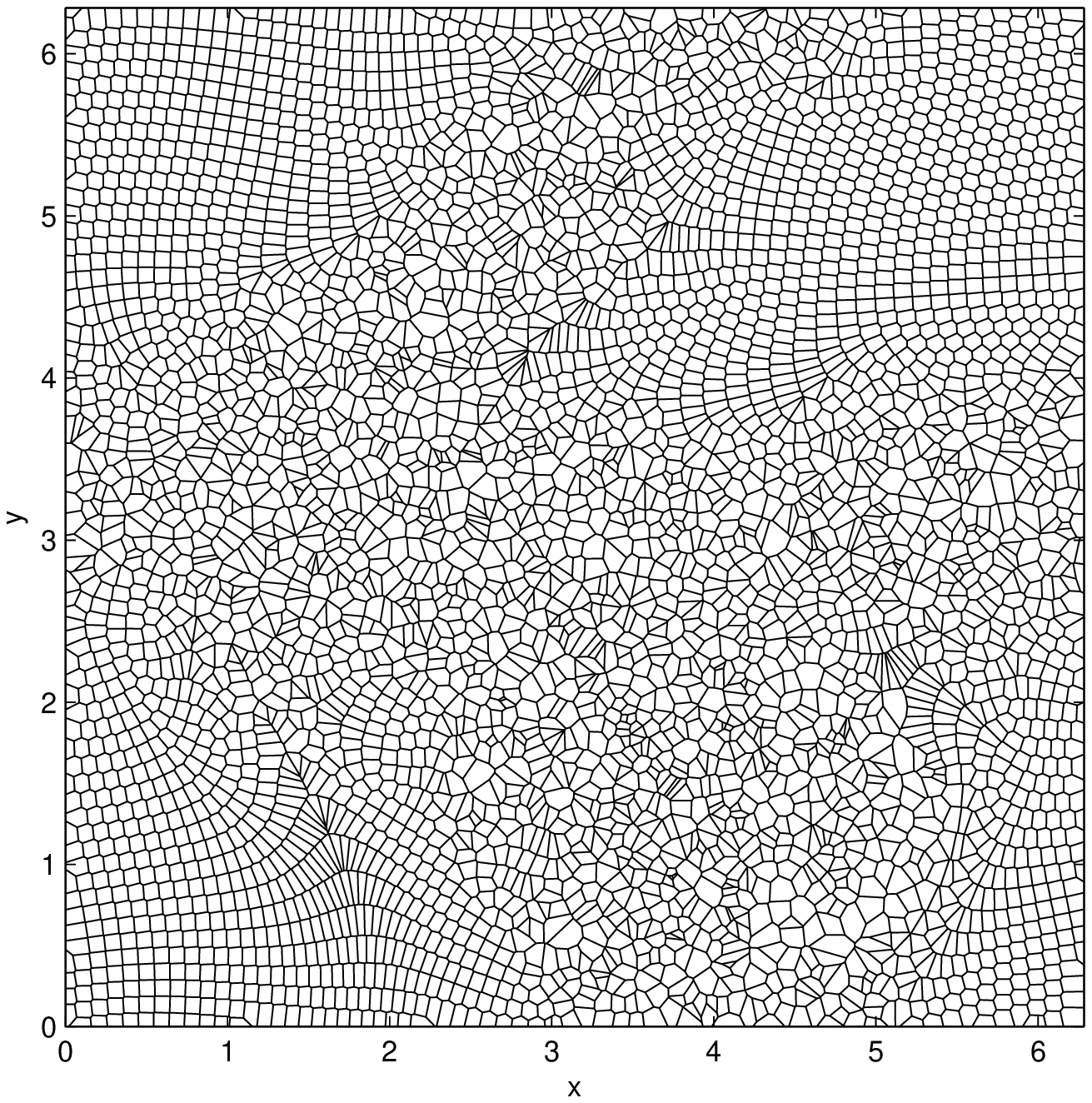}
\includegraphics[angle=0.0,width=0.3\textwidth]{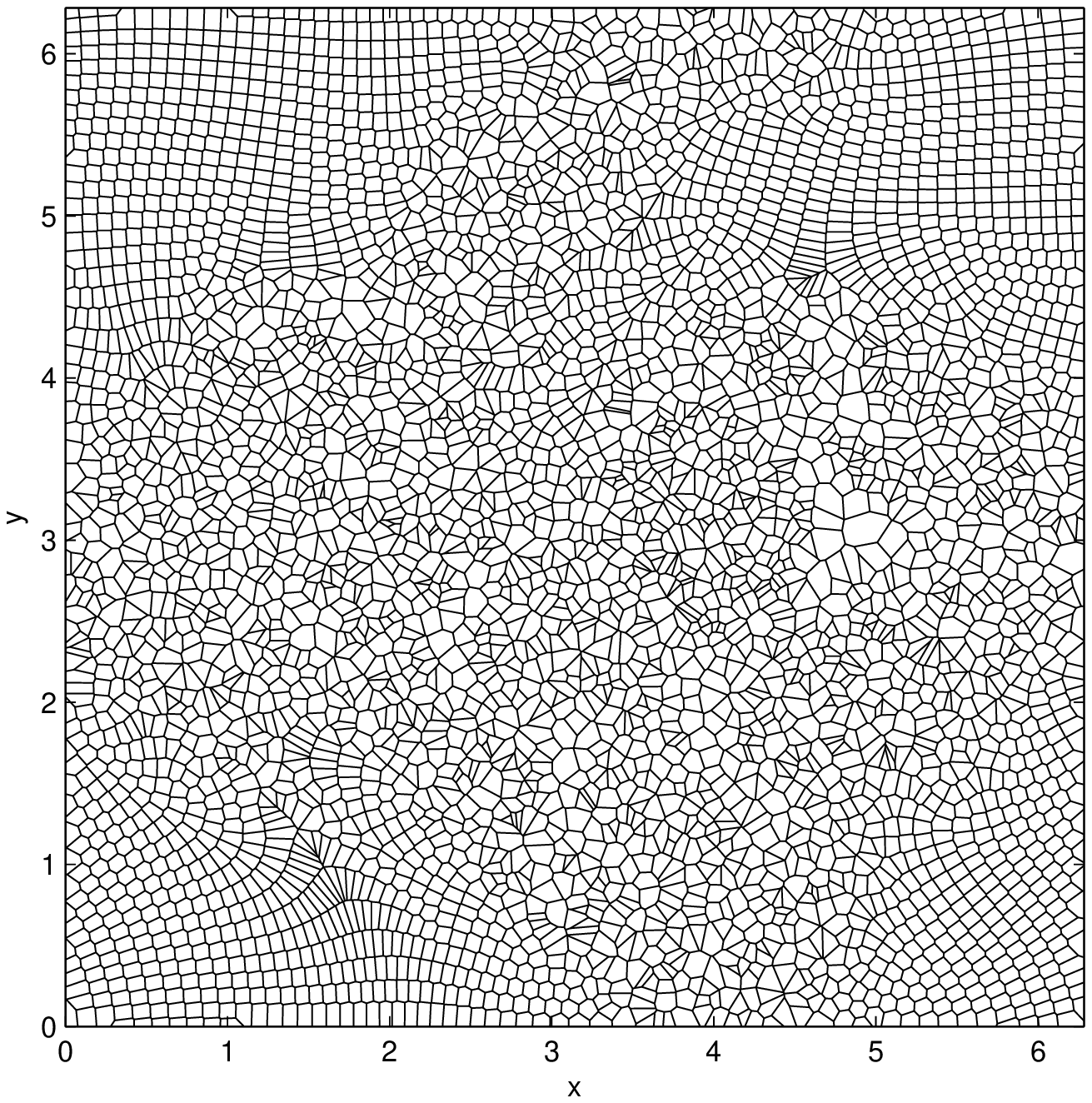}\\
\includegraphics[angle=0.0,width=0.3\textwidth]{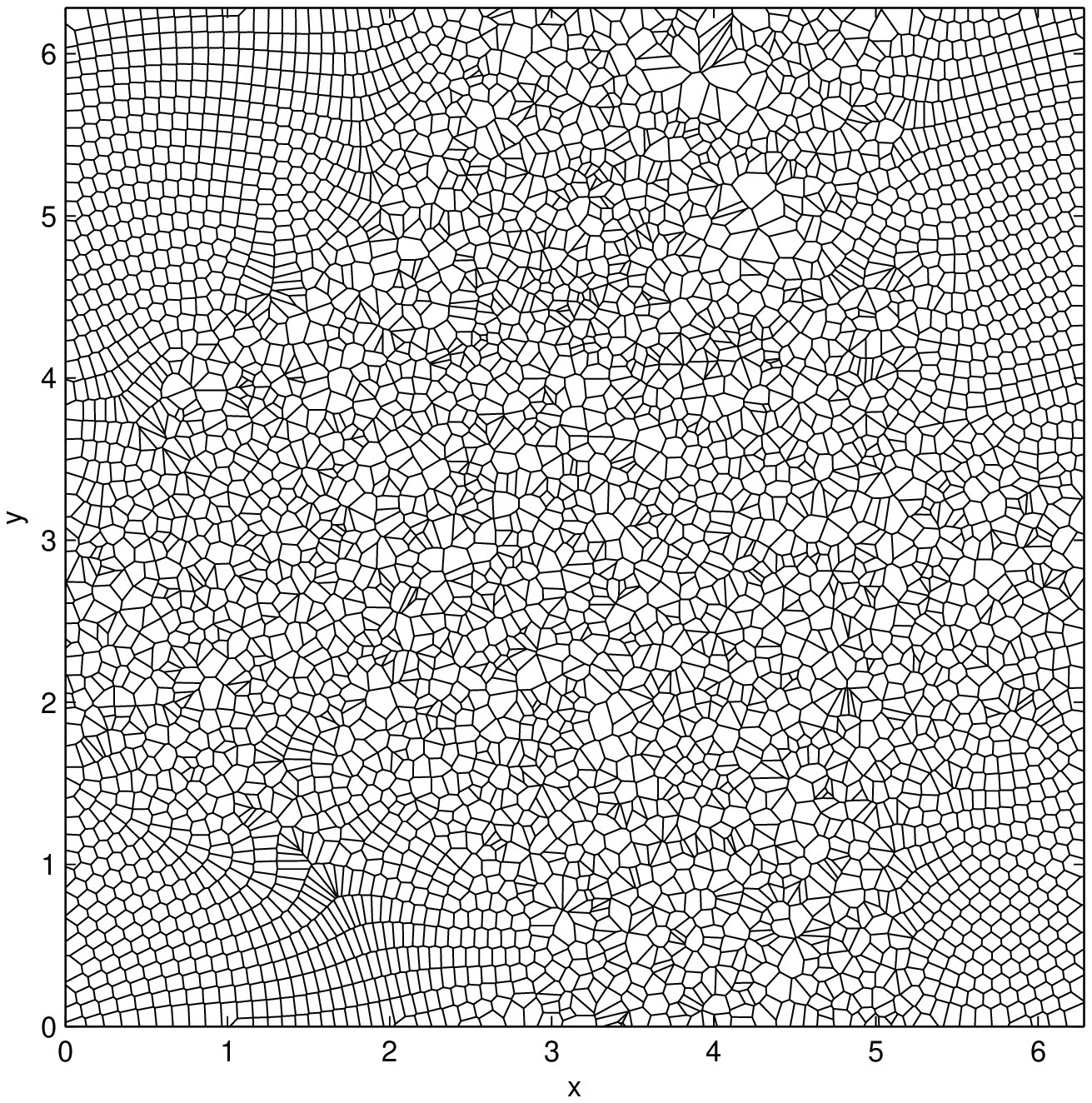}
\includegraphics[angle=0.0,width=0.3\textwidth]{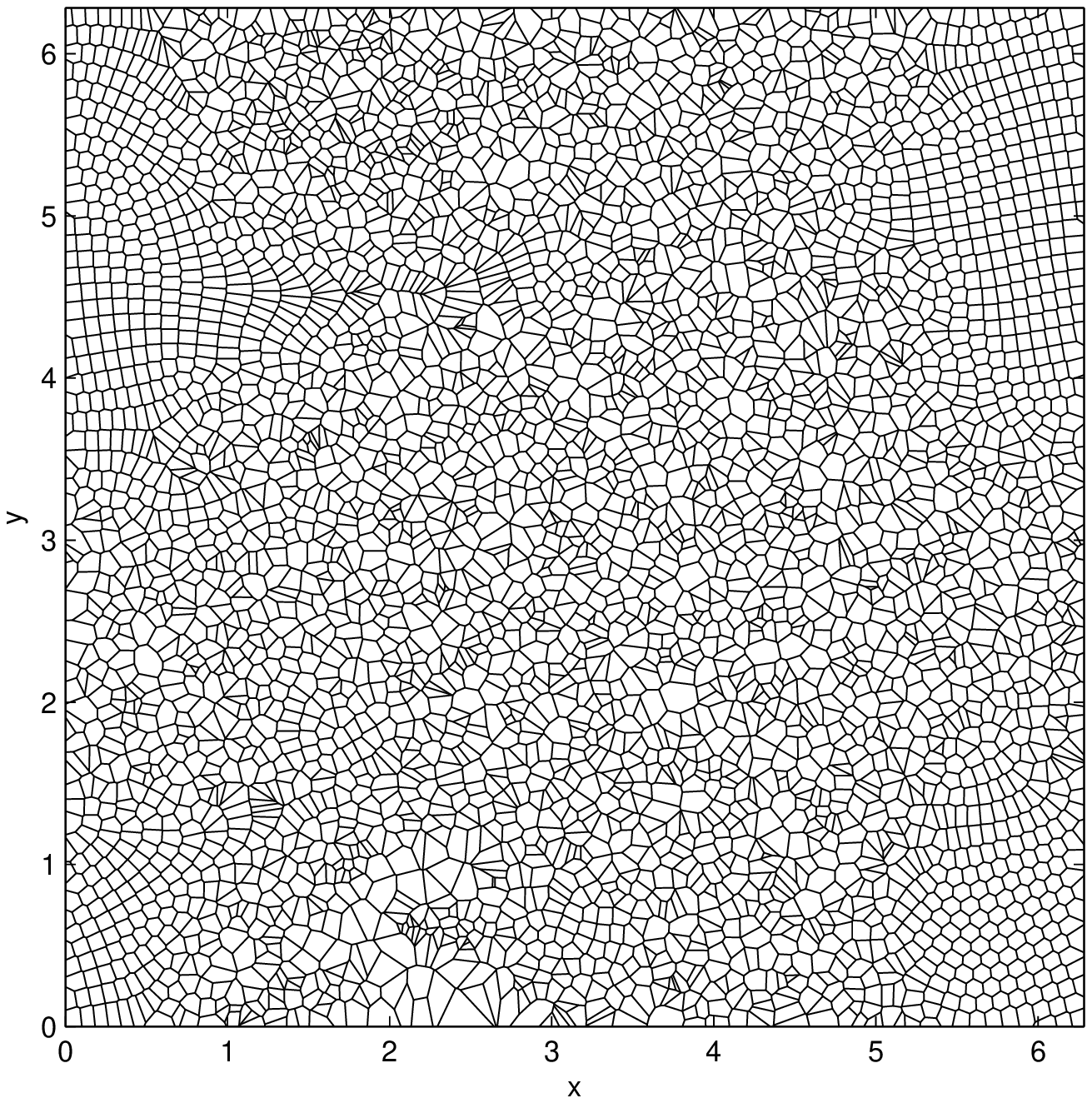}\\
\end{center}
\caption{This graph shows the Voronoi diagram at increments of $40$ days during simulation of a vortex pair in rotating shallow water over a doubly periodic domain. The simulation parameters for this experiment are given in Table \ref{tab:2d_conditions}.}
\label{fig:vfl_vortex_voronoi}
\end{figure}

A contour plot of the regularized layer thickness is shown in Figure \ref{fig:vfl_vortex} in which each snapshot is taken at increments of 40 days and is viewed from left to right, starting at the top of the page. The vortex pair is observed to undergo pure rotation without significant shape deformation or a change in strength of the poles. The Voronoi diagram corresponding to these layer thickness snapshots is shown in Figure \ref{fig:vfl_vortex_voronoi}.  The diagram retains history of the flow field and largely deforms in the wake of the vortex pair to exhibit some intriguiging rotational structure with loss of the radial symmetry. Over the course of the simulation, an increasing number of cells deform from their initial state and the resulting effect is a more homogeneously deformed diagram.

\begin{figure}[!ht] 
\begin{center}
\includegraphics[angle=0.0,width=0.7\textwidth]{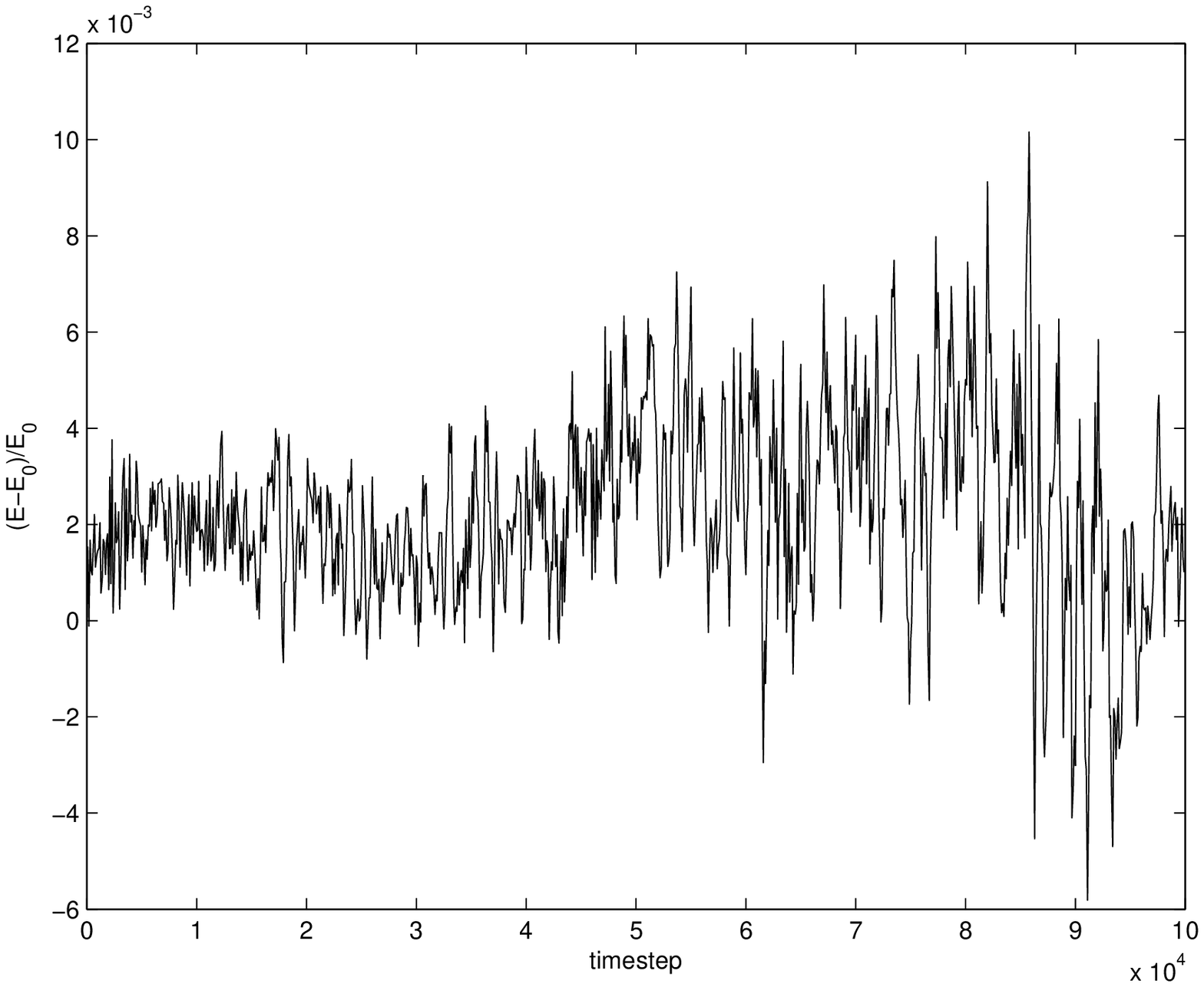}\\
\includegraphics[angle=0.0,width=0.8\textwidth]{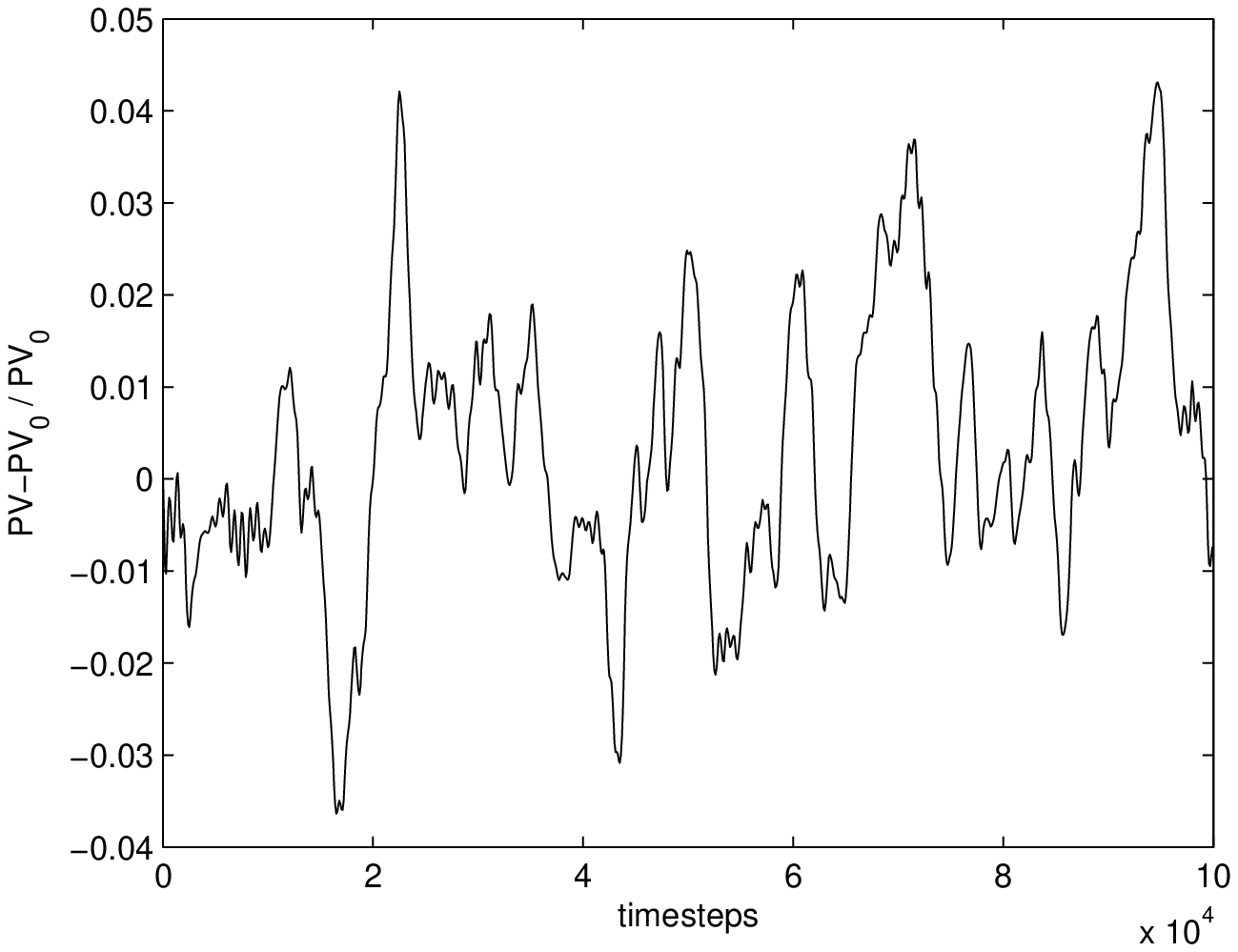}
\end{center}
\caption{(top) The relative energy error and (bottom) the relative potential vorticity error over $1\times10^5$ time steps (approximately 50 years). The simulation parameters for this experiment are given in Table \ref{tab:2d_conditions}.}
\label{fig:vfl_vortex_energy}
\end{figure}

 Figure \ref{fig:vfl_vortex_energy} shows the relative energy and potential vorticity error over a period of $1\times10^5$ time steps or approximately 50 years. The relative energy error does not exhibit secular drift and is an $\mathcal{O}(\Delta t^2)$. The relative potential vorticity error is within 5\% and does not exhibit secular drift either.

%
%

\begin{table}[!h]
\begin{center}
\begin{tabular}{lcccr}
\hline
\vline & Parameter & \vline& value &\vline\\
\hline
\vline &$N_t$& \vline&     $1\times 10^5$ (approx. $50$ years)&\vline\\
\vline &$c_0$& \vline&    $10$&\vline\\
\vline &$f_0$& \vline&    $2\pi$&\vline\\
\vline &$dt$& \vline&     $5\times10^{-2}$&\vline\\
\vline &$L$& \vline&      $2\pi$&\vline\\
\vline &$H$& \vline&      $1$&\vline\\
\vline &$N$& \vline&     $64\times64$&\vline\\
\vline &$n_x$& \vline&     $32$&\vline\\
\vline &$n_y$& \vline&     $32$&\vline\\
\vline &$\hat{\alpha}$& \vline& $4dx$&\vline\\
\hline
\end{tabular}
\caption{Experiment 3: the simulation parameters for the rotating shallow water equations over a doubly-periodic domain initialized by perturbing a geostrophically balanced vortex pair.}
\label{tab:2d_conditions}
\end{center}
\end{table}

\clearpage

\clearpage
\clearpage

\section{Conclusion} \label{sect:sum}
The development of new numerical methods for climate modeling is challenged by the simultaneous need to tractably simulate long-time dynamics and exhibit consistent potential vorticity dynamics.
This paper considers the conservative properties of a variational formulation of the free-Lagrange method for the regularized rotating shallow water equations. The VFL method is based on the free-Lagrange method, a Lagrangian method which remedies the diffusive transport computations of Eulerian methods, by building the transport into the evolution of
the cell positions. Our implementation of the method has the distinguishing feature of using both a Voronoi diagram to formulate a local mass conservation law and the dual Delaunay triangulation to approximate a dispersive regularization of the layer thickness. The semi-discrete regularized shallow water equations preserve symplectic structure and, when integrated with a symplectic time stepper, conserve energy over long-time simulations to the order of the symplectic integrator. 

The discrete variational principle provides a unifying point of conferred mathematical understanding from which to systematically derive new geometric methods. Noether's theorem, for example, informs us that exact potential vorticity conservation is only permissable when the discrete variational principle is invariant under \emph{continuous} particle relabelling. In the absence of this property, we turn to the diagnostic semi-discrete potential vorticity dynamics and choose a discrete curl operator which has the property that it operates on the layer thickness gradient to give a zero vector field. This property guarantees the existence of a semi-discrete potential vorticity law. The semi-discrete potential vorticity and the semi-discrete divergence equations kinematically equate to the respective solenoidal and irrotational components of the velocity field. Together, they augment the description of how the momentum evolves in discrete space. Like the continuum divergence equation (expressed in the Lagrangian frame), the semi-discrete divergence equation exhibits a $\text{div}^2\mathbf{U}$ term which indicates that the flow has a very strong tendency to reach a purely rotational state and therefore has a marked effect upon potential vorticity evolution. 

Numerical results
show (i) the preservation of shape and strength of an initially radially symmetric vortex pair in purely rotational regularized 2D shallow water and (ii) how the Voronoi diagram retains the history of the flow field and (iii) that energy is conserved to $\mathcal{O}(\Delta^2)$ and potential vorticity error to within $5\%$ with no clearly deciperable secular growth over a 50 year period. In the 1D case, we show that the potential vorticity and energy error monotonically decreases with a decreasing time step, the latter as a $\mathcal{O}(\Delta^2)$ scaling law. Simulation of a geostrophic adjustment mechanism provides further evidence that the VFL method exhibits consistent potential vorticity dynamics.


  


\vspace{\baselineskip}

\appendix
\section{Additional Calculations} 

\subsection{Calculation of the semi-discrete vorticity equations}\label{sect:vor_calc}

\be
\begin{split}
&\kh\cdot \text{curl}(\dot{\mathbf{U}}_{\alpha})\\
&=\frac{1}{A_{\alpha}}\sum_i \dot{\mathbf{U}}\abi\cdot \mathbf{d\tau}\abi\\
&=\frac{D}{Dt}\left(\frac{1}{A_{\alpha}}\sum_i\mathbf{U}\abi\cdot \mathbf{d\tau}\abi \right) +\frac{1}{A_{\alpha}^2}
\frac{D}{Dt}(A_{\alpha})\sum_i\mathbf{U}\abi\cdot \mathbf{d\tau}\abi -\frac{1}{A_{\alpha}}\sum_i\mathbf{U}\abi\cdot\frac{D}{Dt}\mathbf{d\tau}\abi\\
&=\frac{D}{Dt}\left(\frac{1}{A_{\alpha}}\sum_i\mathbf{U}\abi\cdot \mathbf{d\tau}\abi \right) +\frac{1}{A_{\alpha}^2}
\frac{D}{Dt}(A_{\alpha})\sum_i\mathbf{U}\abi\cdot \mathbf{d\tau}\abi +\frac{1}{A_{\alpha}}\underbrace{\sum_i\mathbf{U}\abi\cdot\Delta U\abi}_{=0}\\
&=\left(\frac{D}{Dt} + \frac{1}{A_{\alpha}}\frac{D}{Dt}A_{\alpha}\right)\frac{1}{A_{\alpha}}\sum_i\mathbf{U}\abi\cdot \mathbf{d\tau}\abi\\
&=\left(\frac{D}{Dt} + \frac{1}{A_{\alpha}}
\frac{D}{Dt}A_{\alpha}\right)\frac{1}{A_{\alpha}}\sum_i\mathbf{U}\abi\cdot \mathbf{d\tau}\abi\\
&=\left(\frac{D}{Dt} + \text{div}(\mathbf{U}_{\alpha})\right)\frac{1}{A_{\alpha}}\sum_i\mathbf{U}\abi\cdot
\mathbf{d\tau}\abi\\
&=\frac{D}{Dt}\left(\kh\cdot \text{curl}(\mathbf{U}_{\alpha})\right) + \text{div}(\mathbf{U}_{\alpha})
\kh\cdot\text{curl}(\mathbf{U}_{\alpha}).\\
\end{split}
\ee

\vspace{1.0cm}

\begin{remark}
The vanishing under-braced term in the fourth line follows from the definition
of $\Delta U\abi$ and the property that the terms cancel around a closed
loop

\be
\sum_i\mathbf{U}\abi \left(\mathbf{U}_{\alpha}^{\beta_{i+1}}-\mathbf{U}_{\alpha}^{\beta_{i-1}}\right)=0.
\ee
Note that in the continuum limit, the expression takes the form
\be
\frac{1}{2A}\oint_{A}d|\mathbf{U}|^2=0.
\ee
\end{remark}

\subsection{Calculation of the semi-discrete divergence equations} \label{sect:div_calc}
\be
\begin{split}
&\text{div}(\dot{\mathbf{U}}_{\alpha})=\frac{1}{A_{\alpha}}\sum_i \dot{\mathbf{U}}\abi\cdot \mathbf{dn}\abi\\
&=\frac{D}{Dt}\left(\frac{1}{A_{\alpha}}\sum_i\mathbf{U}\abi\cdot \mathbf{dn}\abi \right) +\frac{1}{A_{\alpha}^2}
\frac{D}{Dt}(A_{\alpha})\sum_i\mathbf{U}\abi\cdot \mathbf{dn}\abi -\frac{1}{A_{\alpha}}\sum_i\mathbf{U}\abi\cdot\frac{D}{Dt}\mathbf{dn}\abi\\
&=\frac{D}{Dt}\left(\frac{1}{A_{\alpha}}\sum_i\mathbf{U}\abi\cdot \mathbf{dn}\abi \right) +\frac{1}{A_{\alpha}^2}
\frac{D}{Dt}(A_{\alpha})\sum_i\mathbf{U}\abi\cdot \mathbf{dn}\abi -\frac{1}{A_{\alpha}}\underbrace{\sum_i\mathbf{U}\abi\cdot\Delta(\mathbf{U}^{\perp})\abi}_{:=\Gamma_{\alpha}\neq0}\\
&=\left(\frac{D}{Dt} + \frac{1}{A_{\alpha}}\frac{D}{Dt}A_{\alpha}\right)\frac{1}{A_{\alpha}}\sum_i\mathbf{U}\abi\cdot \mathbf{dn}\abi -\Gamma_{\alpha}\\
&=\left(\frac{D}{Dt} + \frac{1}{A_{\alpha}}
\frac{D}{Dt}A_{\alpha}\right)\frac{1}{A_{\alpha}}\sum_i\mathbf{U}\abi\cdot \mathbf{dn}\abi -\Gamma_{\alpha}\\
&=\left(\frac{D}{Dt} + \text{div}(\mathbf{U}_{\alpha})\right)\frac{1}{A_{\alpha}}\sum_i\mathbf{U}\abi\cdot
\mathbf{dn}\abi - \Gamma_{\alpha}.\\
\end{split}
\ee

\end{document}